\documentclass{article}
\pdfoutput=1
\usepackage{amssymb,amsmath}
\usepackage{amsfonts}
\usepackage{latexsym}
\newtheorem{Lem}{Lemma}
\newtheorem{Thm}{Theorem}
\newtheorem{Pro}{Proposition}
\newtheorem{Rem}{Remark}
\newtheorem{Cor}{Corollary}
\newtheorem{Exa}{Example}
\newtheorem{Def}{Definition}
\def\C{{\mathbb C}}
\def\R{{\mathbb R}}

\def\H{{\mathbb H}}
\def\M{{\mathbb M}} 
\def\S{{\mathbb S}}
\def\z{{\bf z}}

\title{Analytic Functional Calculus in Quaternionic Framework}
\author{Florian-Horia Vasilescu\\
\small Department of Mathematics, University of Lille,\\
\small 59650 Villeneuve d'Ascq, France\\
\small florian.vasilescu@univ-lille.fr}
\date{December 6, 2018}

\begin{document}

\maketitle

\begin{abstract} Regarding quaternions as normal matrices, we first characterize the $2\times 2$ matrix-valued functions, defined on subsets of quaternions, whose values are quaternions. Then we investigate the regularity of quaternionic-valued functions, defined by the analytic functional calculus. 
Constructions of analytic functional calculi for real linear operators, in particular for quaternionic linear ones, are finally discussed, using a Riesz-Dunford-Gelfand type kernel in 
one variable, or a Martinelly type kernel in two variables. 
\end{abstract}
\medskip

{\it Keywords:} quaternionic valued functions; analytic functional calculus;\\ quaternionic spectrum; real operators

{\it Mathematics Subject Classification} 2010: 30G35; 30A05; 47A10; 47A60

\section{Introduction}\label{I}

Introduced in science by W. R. Hamilton as early as 1843,
the quaternions form a unital non commutative division algebra, with numerous applications in mathematics and physics. 
In mathematics, the celebrated Frobenius theorem, proved in 1877, placed the algebra of quaternions among the only three
finite dimensional division algebras over the real numbers,
which is a remarkable feature shared with the real and complex fields.  

Concerning
physics, one can find a first suggestion of a quaternion quantum mechanics in a footnote of a 
paper by   G. Birkhoff and J. von Neumann (see \cite{BiNe}), as mentioned in \cite{FiJaScSp}. 
The work \cite{FiJaScSp} itself presents a quaternionic quantum 
mechanics, using various entities assuming quaternionic values.
The existence of serious physical hypotheses incited several 
mathematicians to develop a branch of analysis in the framework
of quaternions. 

One of the most important investigation in the quaternionic context has been to find a convenient manner to express the 
''analyticity`` of functions depending on quaternions.  Among the pioneer contributions in this direction one should mention
the works \cite{MoTh} and \cite{Fue}. 

More recently, a concept of {\it slice regularity} for functions of one quaternionic variable has been
introduced  in \cite{GeSt}, leading to a large development 
sythesized in \cite{CoSaSt} (which contains a large list of references), whose impact is still actual (see \cite{GhMoPe}, \cite{CoGaKi}, etc.).

Unlike in \cite{GeSt}, the basic idea of the present paper 
is to define the regularity
of a quaternionic-valued function via the analytic functional calculus acting on quaternions. We have chosen to consider the 
algebra of quaternions not as an abstract object but as a
real subalgebra of the complex algebra of $2\times 2$ matrices with complex entries. This (classical) representation has been 
already used by the present author in \cite{Vas4,Vas5}, and 
appears in many other works. 
 Among the advantages of this representation is that we may view  the quaternions as linear operators actually on complex spaces, commuting with the complex numbers. Another one is to regard each quaternion as a normal operator, having a spectrum which can be used to define various compatible functional calculi, including the analytic one. A suggestive parallel is the study of the abstract $C^*$-algebras
as subalgebras of bounded operators on some Hilbert spaces.

One of the main results of this work is Theorem \ref{sym_spec}, giving a characterization of those matrix-valued functions,
defined on some open sets in the complex plane, producing quaternions when applied, by functional calculus, to quaternions having spectra in their domain of definition. Such a  function, temporarily called  {\it skew conjugate symmetric},  corresponds to the more known concept of {\it stem function} (a notion going back  to \cite{Fue}), transposed in our framework  (see Remarks \ref{ex_cex} and \ref{func_scs}).

 Roughly speaking, and unlike in \cite{GeSt}, a ''quaternionic regular function`` can and will be obtained by a pointwise application of  the analytic functional calculus with stem functions on a conjugate symmetric open set $U$ in the complex plane, to quaternions whose spectra are in $U$,  via the  matrix version of  Cauchy's formula (\ref{Cauchy_vect}), with no need of slice derivatives.
In this way, we obtain a whole class of ''regular functions`` (in fact, quaternionic Cauchy transforms of stem functions), having  some unexpected multiplicative properties (see Theorem \ref{H_afc}). In addition, we can and will  recapture, as an illustration,  with our methods and in our terms, several properties of slice regular functions (see Lemma \ref{Cauchy_ineg}, Remark \ref{zeros} etc.), presented for example in  \cite{CoSaSt}. 

The discussion concerning the ''regularity`` of the quaternionic-valued functions is ended with a comaparison between our concept of regularity with  that of slice regularity (Theorem \ref{equiv-ons-dom}), showing that these concepts are equivalent on open sets called in this work {\it spectrally saturated}, which happen to be axially symmetric sets, introduced in \cite{CoSaSt} (see Proposition \ref{UH_tildeU}).

One of the necessities of the quantum mechanics is the existence
of a convenient operator theory in the context of quaternions,
in particular an appropriate definition of the spectrum.  
Because the direct extension of the classical definition of the spectrum has been considered not to be satisfactory, a different definition, using the squares of operators, has been
%parasite text:  bewhose associated functional calculus, when %regarded as normal  matrices, produces quaternions.en 
introduced in \cite{CoSaSt} (see also \cite{CoGaKi}).
We adapt this definition to our framework (see Definition
\ref{Q-spectrum}), showing that, nevertheless,
the symetrization of the classical definition of the 
spectrum is still usable. In fact, we consider a spectrum for real operators on real Banach spaces, and  construct an analytic
functional calculus for them (see Theorem \ref{afcro1} and
Proposition \ref{afcro2}). Unlike in  \cite{CoSaSt} or \cite{CoGaKi}, our functional calculus is based on a Riesz-
Dunford-Gelfand formula, defined in a commututative context,
rather than the non-commutative Cauchy type formula, used by 
previous authors. This analytic functional calculus 
holds for a class of analytic operator valued functions, whose definition extends that of stem functions (see Remark \ref{new_hol_sp}), and it applies, in particular, to some quaternionic linear operators (Corollary \ref{afcqo}).

 \section{Hamilton's Algebra}\label{HA}

We start this discussion with some well known facts. Abstract
Hamilton's algebra $\H_0$ is the four-dimensional $\R$-algebra with 
unit $1$, generated by $\{\bf{j,k,l}\}$, where $\bf j,k,l$ 
satisfy
$$
{\bf jk=-kj=l,\,kl=-lk=j,\,lj=-jl=k,\,jj=kk=ll}=-1.
$$

In this work, {\it Hamiltonn's algebra}  (or {\it algebra of quaternions}) will be identified with an $\R$-subalgebra of 
$\M_2$ of $2\times 2$ matrices with complex
entries. Specifically, using a well-known idea, one considers
the following $2\times2$-matrices with complex entries 

\begin{displaymath}
\bf{I}=
\left( \begin{array}{cc} 1 & 0 \\ 0 & 1\end{array}\right),\quad
\bf{J}=
\left( \begin{array}{cc} i & 0 \\ 0 & -i\end{array}\right),\quad
\bf{K}=
\left( \begin{array}{cc} 0 & 1 \\ -1 & 0\end{array}\right),\quad
\bf{L}=
\left( \begin{array}{cc} 0 & i \\ i & 0\end{array}\right),
\end{displaymath}
with $i^2=-1$. Because we have 
$$
{\bf J}^2={\bf K}^2={\bf L}^2=-{\bf I},$$
$$\bf JK=L=-KJ,\,KL=J=-LK,\,LJ=K=-JL,$$
the assignment
\begin{equation}
\H_0\ni x_0+x_1{\bf j}+x_2{\bf k}+x_3{\bf l}\mapsto x_0{\bf I}+x_1{\bf J}+x_2{\bf K}+x_3{\bf L}\in\M_2
\end{equation}
is an injective unital $\R$-algebra morphism, which is also an 
isometry. For this reason, from now on, the  algebra of quaternions, denoted by $\H$, is defined as the 
$\mathbb R$-subalgebra of the algebra ${\M}_2$, generated by the matrices $\bf I$, $\bf J$, $\bf K$ and $\bf L$. Although this 
identification is not canonical, the realization
of $\H_0$ as a matrix algebra $\H$ offers more freedom when acting with its elements, as we shall see in the sequel. In particular, the quaternions in this paper  commute with complex scalars because they are, in fact, {\it matricial quaternions}. 

We shall often use the notation
\begin{displaymath}
Q({\bf z})=
\left( \begin{array}{cc} z_1 & z_2 \\ -\bar z_2 & \bar z_1\end{array}\right)
\end{displaymath}
for every ${\bf z}=(z_1,z_2)\in\mathbb{C}^2$, and noticing that
$$Q({\bf z})=\Re{z}_1{\bf I}+\Im{z}_1{\bf J}+\Re{z}_2{\bf K}+ \Im{z}_2{\bf L},$$
and ${\bf I}=Q((1,0)),\,{\bf J}=Q((i,0)),\,{\bf K}=Q((0,1)),\,{\bf L}=Q((0,i))$, we obtain that the map 
$\C^2\ni{\bf z}\mapsto Q({\bf z})\in\H$ is $\R$-linear
and bijective. In other words, the set $\C^2$ can be identified, as an $\R$-vector space, with the algebra $\H$. For technical reasons,
we often represent a fixed element of $\H$ under the form 
$Q(\z)$, for some $\z\in\C^2$ uniquely determined, via the assignment (1).  

Regarding the elements of $\M_2$ as linear maps acting on the 
space $\C^2$, endowed with the natural scalar product 
$\langle{\bf z,w}\rangle=z_1\bar{w}_1+z_2\bar{w}_2$ and the 
associated norm 
$\Vert{\bf z}\Vert^2=\vert z_1\vert^2+\vert z_2\vert^2,\, {\bf z}=(z_1,z_2),{\bf w}=(w_1,w_2)\in\C^2$,  we see that
the algebra $\H$ also has a natural involution, given by 
$Q(\z)\mapsto Q(\z)^*,\,\z\in\C^2,$ where
%%%%%%%%%%%%%%%%%%%%%%%%%%%%%%%%%%%%%%%%%%%%%%%%%%%%%%%%%%%%%%
% Text parasite
%, for , 
%which is the same as saying that $F$ is a stem function
%%%%%%%%%%%%%%%%%%%%%%%%%%%%%%%%%%%%%%%%%%%%%%%%%%%%%%%%%%%%
${\bf z}=(z_1,z_2)$,
\begin{displaymath}
Q({\bf z})^*=
\left( \begin{array}{cc} \bar{z}_1 & -z_2 \\ \bar z_2 &  z_1\end{array}\right)=Q({\bf z}^*),
\end{displaymath}
with ${\bf z}^*=(\bar{z}_1,-z_2)$. (In fact, the map $\C^2\ni\z\mapsto\z^*\in\C^2$ is itself an $\R$-linear involution of $\C^2$.) 
In particular, ${\bf J}^*=-{\bf J},\,{\bf K}^*=-{\bf K},\,{\bf L}^*=-{\bf L}.$
 
It is easily seen  that $Q({\bf z})Q({\bf z})^*=Q({\bf z})^*Q({\bf z})=\Vert {\bf z}\Vert^2{\bf I}$ for all 
${\bf z}\in{\mathbb C}^2$, and so $Q({\bf z})$ is a normal matrix for each ${\bf z}\in{\mathbb C}^2$.
Moreover,  $\Vert Q({\bf z})\Vert=\Vert {\bf z}\Vert$ for all ${\bf z}\in{\mathbb C}^2$, that is, the map  $\C^2\ni\z\mapsto Q(\z)\in\H$ is an isometry. In addition,  
$Q({\bf z})^{-1}=\Vert {\bf z}\Vert^{-2}Q({\bf z})^*$ for all 
${\bf z}\in{\mathbb C}^2\setminus\{0\}$, so every nonnull element of $\H$ is invertible. 

Let us clarify the position of the $\R$-subalgebra $\H$ into the  $\C$-algebra $\M_2$.

% Rem 1

\begin{Rem}\label{skewcc}\rm On the algebra $\M_2$ we define what we will call a  {\it skew complex conjugation}, setting
\begin{displaymath}
{\bf a}^\sim:=
\left( \begin{array}{cc} \bar{a}_4 & -\bar{a}_3 \\ -\bar{a}_2 & \bar{a}_1\end{array}\right),
\end{displaymath}
for every
\begin{displaymath}
{\bf a}=
\left( \begin{array}{cc} a_1 & a_2 \\ a_3 & a_4\end{array}\right)
\in\M_2.
\end{displaymath}
The map ${\bf a}\mapsto{\bf a}^\sim$ is conjugate homogeneous
and additive, in particular  $\mathbb R$-linear, multiplicative, unital, $({\bf a}^\sim)^\sim={\bf a}$,  and $({\bf a}^*)^\sim=({\bf a}^\sim)^*$. In addition, ${\bf a}={\bf a}^\sim$ 
% Misprint: $\tilde{\bf a}$ replaced by ${\bf a}^\sim$
if and only if ${\bf a}$ is a quaternion.
 
Being a $*$-automorphism $\R$-linear, the map ${\bf a}\mapsto{\bf a}^\sim$ must be an isometry. Because $\M_2$ is finite dimensional, we can give a direct easy argument of this assertion, as follows. Because
$$
{\bf a}^*{\bf a}=\left(\begin{array}{cc} \vert a_1\vert^2+\vert a_3\vert^2 &
\bar{a}_1a_2+\bar{a}_3a_4\\ \bar{a}_2a_1+\bar{a}_4a_3 & 
\vert a_2\vert^2+\vert a_4\vert^2\end{array}\right),
$$
$$
({\bf a}^\sim)^*{\bf a}^\sim=\left(\begin{array}{cc}\vert a_2\vert^2+\vert a_4\vert^2 & -\bar{a}_1a_2-\bar{a}_3a_4 \\ -\bar{a}_2a_1-
\bar{a}_4a_3 & \vert a_1\vert^2+\vert a_3\vert^2\end{array}\right),
$$
we derive the equality
${\rm det}(\lambda{\bf I}-({\bf a}^\sim)^*{\bf a}^\sim)={\rm det}(\lambda{\bf I}-{\bf a}^*{\bf a}),$  for all $\lambda\in\C$,
so the matrices ${\bf a}^*{\bf a},\, ({\bf a}^\sim)^*{\bf a}^\sim$ have the same spectrum.
Because the norms of the positive $({\bf a}^\sim)^*{\bf a}^\sim,\, {\bf a}^*{\bf a}$ equal the greatest joint eigenvalue, we should
have $\Vert{\bf a}^\sim\Vert^2=\Vert ({\bf a}^\sim)^*{\bf a}^\sim\Vert=\Vert {\bf a}^*{\bf a}\Vert=\Vert {\bf a}\Vert^2$. 
\end{Rem}

Note also that 
$$
{\bf a}=\frac{{\bf a}+{\bf a}^\sim}{2}
+i\frac{{\bf a}-{\bf a}^\sim}{2i},\,\,{\bf a}\in\M_2,
$$
with ${\bf a}+{\bf a}^\sim,\,i({\bf a}-{\bf a}^\sim)\in\H$. 
In other words, $\M_2=\H+i\H$.  We also have $\H\cap i\H=\{0\}$. 
Indeed, if $q=ir$ with $q,r\in\H$, we have $q^\sim=q=
(ir)^\sim=-ir=-q$, whence $q=0$, showing that the decomposition
$\M_2=\H+i\H$ is a direct sum.

A map similar to the skew complex conjugation is defined, under 
the name of {\it reflexion}, in $C^*$-algebras (see \cite{LorSor}, Definition 2.6.). Nevertheless, a reflexion is an anti-automorphism, by definition.

\section{A Spectral Approach to $\H$-Valued Functions}
\label{SAH-VF}

 As before, the space $\C^2$ is endowed with its natural scalar product $\langle*,*\rangle$, and norm $\Vert*\Vert$. We also have the $\R$-linear the map 
 $\C^2\ni{\bf z}\mapsto Q({\bf z})\in\H,$ which 
% added: is a 
is a  bijective isometry. In other words, giving a quaternion $q\in\H$, there is a unique point $\z_q\in\C^2$ such that $q=Q(\z_q)$. Moreover, the algebra $\H$ will be regarded as an $\R$-subalgebra of the $\C$-algebra $\mathbb{M}_2$. In this way, the elements of
$\mathbb{H}$ will be regarded as linear operators on the 
Hilbert space $\C^2$, so every element $Q({\bf z})$ is a normal
operator on the Hilbert space $\C^2$. 

Occasionally, we use the notation $\Re q=\Re z_1$ and $\Vert q\Vert=
\Vert Q(\z)\Vert$ if $q=Q(\z)$ and $\z=(z_1,z_2)$. 

Let us mention that for every complex space Banach 
$\mathcal{X}$, and each Banach space operator $T$ on 
$\mathcal{X}$, in this text the symbol $\sigma(T)$ will designate the spectrum of $T$, and the symbol $\rho(T)$ will be
resolvent set of $T$. Other similar symbols will be later
introduced and explained. 

We start with an elementary result:

% Lem 1

\begin{Lem}\label{roots} Let $\z=(z_1,z_2)\in\C^2$ be fixed. The spectrum
$\sigma(Q({\bf z}))=\{s_{\pm}({\bf z})\}$ of the normal operator $Q({\bf z})$  is given by
\begin{equation}\label{spectrum}
s_{\pm}({\bf z})=\Re z_1\pm i\sqrt{(\Im z_1)^2+\vert z_2\vert^2},\,\,
{\bf z}=(z_1,z_2)\in\C^2.
\end{equation}
We have $s_+(\z)=\overline{s_-(\z)}$, and the  points $s_+(\z),s_-(\z)$ are distinct if and only if 
$Q(\z)\notin \R{\bf I}$. Moreover:

{\bf (a)}   if $z_2\neq0$, the elements
\begin{equation}\label{val_prop}
\nu_{\pm}(\z)=\frac{1}{\sqrt{\vert z_2\vert^2+\vert s_{\pm}(\z)-z_1\vert^2}}(z_2,s_{\pm}(\z)-z_1)\in\C^2,\,\,\z=(z_1,z_2)\in\C^2
\end{equation}
are 
% "the" has been erased 
eigenvectors corresponding to the eigenvalues $\{s_{\pm}({\bf z})\}$ respectively, and they
form an orthonormal basis of the Hilbert space $\C^2$; 

{\bf (b)} if $z_2=0$ but $\Im z_1\neq 0$, we have $\sigma(Q(\z))=\{s_{\pm}({\bf z})\}$, with $s_+(\z)=z_1,s_-(\z)=\bar{z}_1$, and $\nu_+(\z)=(1,0),\,\nu_-(\z)=(0,1)$ are eigenvectors corresponding the eigenvalues $z_1,\bar{z}_1$, respectively;

{\bf (c)} if $z_2=0$ and  $\z=(x,0)$ with 
$x\in\R$, we have $\sigma(Q({\bf z}))=\{x\}$, with
$s_+(\z)=s_-(\z)=x$, and $\nu_+(\z)=(1,0),\,\nu_-(\z)=(0,1)$ are eigenvectors corresponding to the eigenvalue $x$. 
\end{Lem}

{\it Proof.}\, The spectrum $\sigma(Q({\bf z}))$ of 
 $Q({\bf z})$ is  given by the roots of the equation
\begin{equation}\label{eq_spectrum}
s^2-2s\Re z_1+\vert z_1\vert^2+\vert z_2\vert^2=0,
\end{equation}
leading to the equality (\ref{spectrum}). 

{\bf (a)} Let  ${\bf z}=(z_1,z_2)\in\C^2$ be with  $z_2\neq0$, so  
$Q(\z)\notin \R{\bf I}$. In this case, clearly
$\vert z_2\vert^2+\vert s_{\pm}(\z)-z_1\vert^2>0$.
The vectors $(z_2,s_{+}(\z)-z_1),\,(z_2,s_{-}(\z)-z_1)$ are 
orthogonal eigenvectors of $Q(\bf z)$ in $\C^2$, corresponding to the eigenvalues $\nu_+(\z),\nu_-(\z)$,  
via equation (\ref{eq_spectrum}). Hence $\{\nu_+(\z),\nu_-(\z)\}$ is an orthonormal basis of $\C^2$, with the stated properties. 

The assertions {\bf (b)}, {\bf (c)} are easily obtained and left 
to the reader. We only note that $z_2=0$ implies, in fact,
$s_{\pm}({\bf z})=\Re z_1\pm i\vert\Im z_1\vert$, leading to 
the eigenvectors $(\Re z_1+i\vert\Im z_1\vert,0)$ and 
$(0,\Re z_1-i\vert\Im z_1\vert)$, corresponding to the  eigenvalues $z_1,\bar{z}_1$, respectively. For the sake of 
simplicity, we take $s_+(\z)=z_1,\,s_-(\z)=\bar{z}_1$, and
replace the eigenvectors from above by $(1,0)$ and $(0,1)$, respectively, with no loss of generality. 

\begin{Exa}\label{spec_can_rep}\rm Let $\S=\{\mathfrak{s}=
x_1{\bf J}+
x_2{\bf K}+x_3{\bf L};x_1,x_2,x_3\in\R,x_1^2+x_2^2+x_3^2=1\}$, that is, the unit sphere of purely imaginary quaternions. As noticed in in \cite{CoSaSt} (and easily seen), every quaternion $q\in \H\setminus\R$ 
can be  written as $q = x{\bf I} + y\mathfrak{s}$,  for 
% The word "each" is replaced by "some".
 some $\mathfrak{s}\in\S$, where $x,y$ are real numbers, which are unique when $y > 0$. Let us prove that for every 
$q=x{\bf I}+y\mathfrak{s},\,x,y\in\R,$ we have $\sigma(q)=\{x\pm iy\}$. 

The quaternion $\mathfrak{s}\in\mathbb{S}$ can be written as a matrix under the form
$$
\mathfrak{s}=\left(\begin{array}{cc} ia_1 & a_2+ia_3\\ -a_2+ia_3
& -ia_1\end{array}\right),
$$
where $a_1,a_2,a_3$ are real numbers with $a_1^2+a_2^2+a_3^2=1$.
Hence, 
$$
\lambda{\bf I}-x-y\mathfrak{s}=\left(\begin{array}{cc} 
\lambda-x-ia_1y & -(a_2+ia_3)y\\ (a_2-ia_3)y & \lambda-x+ia_1y\end{array}\right).
$$
The detrminant of the matrix from above 
% one "is" has been erased
is equal to  $\lambda^2-2\lambda x+x^2+y^2$,
via the equality $a_1^2+\vert a_2+ia_3\vert^2=1$, whose roots are equal to $x\pm iy$. 

Note that the spectrum of $q$ does not depend on $\mathfrak{s}$.

\end{Exa}

\begin{Def}\label{canon_ev}\rm For a fixed point $\z=(z_1,z_2)\in\C^2$, let $\{s_\pm(\z)\}$ be the spectrum of the operator $Q({\bf z})$, given by Lemma \ref{roots}.  The eigenvectors 
$\{\nu_{\pm}(\z)\}$ of $Q({\bf z})$ corresponding to the eigenvalues $\{s_\pm(\z)\}$
respectively, again given by Lemma \ref{roots}, will be called the {\it canonical eigenvectors} of $Q({\bf z})$.

For an arbitrary quaternion $q$, we define its {\it spectrum} 
$\sigma(q)$ as the set equal to $\sigma(Q({\bf z}))=\{s_{\pm}({\bf z})\}$, where ${\bf z}=(z_1,z_2)\in\C^2$ is the unique point with $q=Q({\bf z})$. We also use the notation $s_{\pm}(q)=
s_{\pm}({\bf z})$ and $\nu_{\pm}(q)=\nu_{\pm}({\bf z})$.
\end{Def}

% Lem 2 

\begin{Lem}\label{nu_proper} Let $\z=(z_1,z_2)\in\C^2$, and let $\nu_{\pm}(\z)=(\nu_{\pm1}(\z),\nu_{\pm2}(\z))\in\C^2$ be 
the canonical eigenvectors of $Q(\z)$. Then we have 
\begin{displaymath}\begin{array}{cc}
\vert\nu_{-1}(\z)\vert^2=\vert\nu_{+2}(\z)\vert^2,\,
\vert\nu_{-2}(\z)\vert^2=\vert\nu_{+1}(\z)\vert^2, \\
\hskip10cm(*)\\
\nu_{-1}(\z)\overline{\nu_{-2}}(\z)+\nu_{+1}(\z)\overline{\nu_{+2}}(\z)=0
\end{array}
\end{displaymath}
\end{Lem}

{\it Proof.} We set $\nu_\pm:=\nu_{\pm}(\z)$
and $\nu_\pm=(\nu_{\pm1},\nu_{\pm2})\in\C^2$. If $z_2=0$, we have $\nu_{+1}=\nu_{-2}=1,\,\nu_{-1}=\nu_{+2}=0$, so relations
$(*)$ are trivial.

Assume $z_2\neq0$. Because $\vert s_+-z_1\vert=\vert s_--\bar{z}_1\vert$, and relation (\ref{eq_spectrum}) is equivalent to 
$$
(s_\pm-z_1)(s_\pm-\bar{z}_1)+\vert z_2\vert^2=0, 
$$
we have $\vert s_- -z_1\vert^2\vert s_+-z_1\vert^2=\vert z_2\vert^4$. Therefore
$$\vert\nu_{-1}\vert^2=\frac{\vert z_2\vert^2}{\vert z_2\vert^2+
\vert s_--z_1\vert^2}=\frac{\vert s_+ -z_1\vert^2}{\vert z_2\vert^2+\vert s_+-z_1\vert^2}=\vert\nu_{+2}\vert^2.
$$
A similar argument shows that
$$
\vert\nu_{+1}\vert^2=\frac{\vert z_2\vert^2}{\vert z_2\vert^2+
\vert s_+-z_1\vert^2}=\frac{\vert s_- -z_1\vert^2}{\vert z_2\vert^2+\vert s_--z_1\vert^2}=\vert\nu_{-2}\vert^2.
$$

We also have
$$
\nu_{-1}\overline{\nu_{-2}}=\frac{z_2(s_+-\bar{z}_1)}{\vert z_2\vert^2+\vert s_--z_1\vert^2}=
\frac{z_2(s_+-\bar{z}_1)\vert s_+-z_1\vert^2}{\vert z_2\vert^2
(\vert z_2\vert^2+\vert s_+-z_1\vert^2)}=
$$
$$
-\frac{z_2(s_--\bar{z}_1)}{\vert z_2\vert^2+\vert s_+-z_1\vert^2}=-\nu_{+1}\overline{\nu_{+2}},
$$
via equation (\ref{eq_spectrum}). Consequently, equalities 
$(*)$ hold true.
\medskip

\noindent {\bf Remark} We note that equalities $(*)$ do not follow, in general, from  the orthogonality of $\nu_+(\z)$ and $\nu_-(\z)$.
\medskip

% Rem 2

\begin{Rem}\label{samesp}\rm 
Given a complex number $\zeta$, we can
determine all quaternions $q$ with $\sigma(q)=\{\zeta,\bar{\zeta}\}$. Assuming, with no loss of generality, that $\Im\zeta\ge0$, we look for the points $\z=(z_1,z_2)\in\C$ satisfying the equation
$$
\zeta=s_+(\z)=\Re z_1+i\sqrt{(\Im z_1)^2+\vert z_2\vert^2},
$$
so 
$\bar{\zeta}=s_-(\z)=\Re z_1-i\sqrt{(\Im z_1)^2+\vert z_2\vert^2}.$
Setting $u=z_2$ as a parameter, we obtain $\Re z_1=\Re\zeta$, and 
$(\Im z_1)^2=(\Im\zeta)^2-\vert u\vert^2$, provided $\vert u\vert^2\le (\Im\zeta)^2$. 
The solutions are given by the set
$$
\{\z=(\Re\zeta\pm i\sqrt{(\Im\zeta)^2-\vert u\vert^2},u)\in\C^2,
\,\vert u\vert\le\Im\zeta\},
$$
so we have, for every such a $\z$,   $\sigma(Q(\z))=\{\zeta,\bar{\zeta}\}$, via 
Lemma \ref{roots}. In particular, if $\z=(\Re\zeta + i\sqrt{(\Im\zeta)^2-\vert u\vert^2},u)$ for some $u\in\C$ with 
$0\neq\vert u\vert\le\Im\zeta$, we must have
$$
\left(\begin{array}{cc} \Re\zeta+i\sqrt{(\Im\zeta)^2-\vert u\vert^2} & u\\ -\bar{u} & \Re\zeta-i\sqrt{(\Im\zeta)^2-\vert u\vert^2}\end{array}\right)\left(\begin{array}{cc} u \\
i(\Im\zeta-\sqrt{(\Im\zeta)^2-\vert u\vert^2})\end{array}\right)
$$
$$
=\zeta\left(\begin{array}{cc} u \\
i(\Im\zeta-\sqrt{(\Im\zeta)^2-\vert u\vert^2})\end{array}\right),
$$
which is an explicit form, modulo a multiplicative constant, of the equation $Q(\z)\nu_+(\z)=
s_+(\z)\nu_+(\z)$, with $s_+(\z)=\zeta$. In fact,
$$
\nu_+(\z)=(2\Im\zeta(\Im\zeta-\sqrt{(\Im\zeta)^2-\vert u\vert^2}))^{-1/2}
(u,i(\Im\zeta-\sqrt{(\Im\zeta)^2-\vert u\vert^2})).
$$

If $\Im\zeta\le0$, we apply the previous discussion to $\bar{\zeta}$.

\end{Rem}

% Rem 3

\begin{Rem}\label{consym}\rm (1) A subset $U\subset\C$ is said to be {\it conjugate symmetric} if $\zeta\in U$ if and only if 
$\bar{\zeta}\in U$. 

For an arbitrary conjugate symmetric subset
$U\subset\C$ we put
$$U_\H=\{q\in\H;\sigma(q)\subset U\}.$$ 
Note that, for every $\zeta\in U$ and $u\in\C$ with $\vert u\vert\le\vert\Im\zeta\vert$, setting
$$
q_\zeta^\pm(u):=\Re\zeta\pm i\sqrt{(\Im\zeta)^2-\vert u\vert^2},u)\in\C^2,
\,\vert u\vert\le\vert\Im\zeta\vert,
$$
we have
$$
U_\H=\{Q(q_\zeta^\pm(u));\zeta\in U,\,u\in\C,\,\vert u\vert\le\vert\Im\zeta\vert\},
$$
via Remark \ref{samesp}.

If $U\subset\C$ is open and conjugate symmetric,
the set $U_\H$ is also open via the upper semi-continuity of the spectrum (see \cite{DuSc},
Lemma VII.6.3.). 

(2) A subset $A\subset\H$ is said to be {\it spectrally saturated} 
if whenever $\sigma(r)=\sigma(q)$ for some $r\in\H$ and $q\in A$, 
we also have $r\in A$. 

For an arbitrary $A\subset\H$, we put $\mathfrak{S}(A)=\{\zeta\in\C;
\exists q\in A: \zeta\in\sigma(q)\}$.
As above, we also  put $S_\H=\{q\in\H;\sigma(q)\subset S\}$ for an arbitrary subset $S\subset\C$.

(3) A subset $A\subset\H$ is spectrally saturated if and only if there exists a conjugate symmetric subset $S\subset\C$ such that 
$A=S_\H$. In this case, $S=\mathfrak{S}(A)$. 

Indeed, if $S\subset\C$  is conjugate symmetric and $A=S_\H$, then clearly $A$ is spectrally saturated. Conversely, if $A$ is 
spectrally saturated, then $A=\mathfrak{S}(A)_\H$, because we always 
have $A\subset\mathfrak{S}(A)_\H$, and fixing $r\in\mathfrak{S}(A)_\H$,
we can find $q\in A$ with $\sigma(q)=\sigma(r)$, so $r\in A$. Finally,
if $\mathfrak{S}(A)_\H=S_\H$ for some conjugate symmetric $S\subset\H$,
we must have $S=\mathfrak{S}(A)$. 

(4) If $\Omega\subset\H$ is an open spectraly saturated set, then 
$\mathfrak{S}(\Omega)\subset\C$ is open. Indeed, we have an injective 
$\R$-linear map $\mathfrak{S}(\Omega)\ni\zeta\mapsto Q((\zeta,0))\in\Omega$, and the subset $\{Q((\zeta,0));\zeta\in\mathfrak{S}(\Omega)\}$ is 
open in the vector space $\{Q((\zeta,0));\zeta\in\C\}$. Hence
the subset $U\subset\C$ is open if and only if $U_\H$ is open,
via an assertion in (1), from above.

An important particular case is when $U={\mathbb D}_r:=\{\zeta\in\C;\vert\zeta\vert<r\}$, for some $r>0$. Because the 
norm of the normal operator induced by $q$ on $\C^2$ is equal to its spectral radius, we must have $U_\H=\{q\in\H;
\Vert q\Vert<r\}$.

(5) We finally note that, for a given  conjugate symmetric subset
$U\subset\C$, the set $U_\H$ is precisely the {\it circularization} of $U$, via Proposition \ref{UH_tildeU}, so it is {\it axially symmetric} (see \cite{GhMoPe}, Section 1.1 and \cite{CoSaSt}, Definition 4.3.1).
 Nevertheless, we continue to call such a set spectrally saturated, a name which better reflects our spectral approach.

\end{Rem}

% Rem 4 
 
\begin{Rem}\label{ex_cex}\rm  Let $U\subset\C$ be conjugate symmetric, and let $F:U\mapsto\M_2$. We write 
\begin{displaymath}
F(\zeta)=
\left( \begin{array}{cc} f_{11}(\zeta) & f_{12}(\zeta) \\ f_{21}(\zeta) & f_{22}(\zeta) \end{array}\right),\,\,\zeta\in U,
\end{displaymath}
with $f_{mn}:U\mapsto\C$, $m,n\in\{1,2\}$, and set 
\begin{displaymath}
F^\sim(\zeta)=
\left( \begin{array}{cc} \overline{f_{22}(\zeta)} & -\overline{f_{21}(\zeta)}
\\ -\overline{f_{12}(\zeta)} & \overline{f_{11}(\zeta)} \end{array}\right),\,\,\zeta\in U.
\end{displaymath}
In other words, $F^\sim(\zeta)=(F(\zeta))^\sim$ for all $\zeta\in U$, where  ''${}^\sim$`` designates the skew complex conjugation (see Remark \ref{skewcc}). 

We temporarily say that $F$ is
{\it skew conjugate symmetric} if $F(\bar{\zeta})=F^\sim(\zeta),\,\zeta\in U$. Clearly,
the function $F$ is skew conjugate symmetric if and only if $f_{11}(\zeta)=\overline{f_{22}(\bar{\zeta})}$ and $f_{12}(\zeta)=-\overline{f_{21}(\bar{\zeta})}$, implying
$f_{22}(\zeta)=\overline{f_{11}(\bar{\zeta})}$ and $f_{21}(\zeta)=-\overline{f_{12}(\bar{\zeta})}$ for all 
$\zeta\in U$. In fact, the function $F$ is  skew conjugate symmetric if and only if $F$ has the form
\begin{equation}\label{func_scs}
F(\zeta)=\left(\begin{array}{cc} f_1(\zeta) & f_2(\zeta) 
\\ -\overline{f_2(\bar{\zeta})} & 
\overline{f_1(\bar{\zeta})}  \end{array}\right),\,\,\zeta\in U,
\end{equation}
for some functions $f_1,f_2:U\mapsto\C$.
\end{Rem}

% Rem 5
 
\begin{Rem}\label{ssc_vs_sf}\rm It is interesting to compare the stem functions (see for instance \cite{GhMoPe}, Section 1.1),
with the skew symmetric conjugate functions. To transpose this 
discussion in our context, let us remark that the tensor product $\H\otimes_\R\C$ may be identified
with $\M_2=\H+i\H$, which is a direct sum, via the  isomorphism induced by the decomposition
$$
{\bf a}=\frac{{\bf a}+{\bf a}^\sim}{2}
+i\frac{{\bf a}-{\bf a}^\sim}{2i},\,\,{\bf a}\in\M_2,
$$
with ${\bf a}+{\bf a}^\sim,\,i({\bf a}-{\bf a}^\sim)\in\H$ (see
Remark \ref{skewcc}). The corresponding conjugation of $\M_2$
is in this case ${\bf a}={\bf b}+i{\bf c}\mapsto\bar{\bf a}=
{\bf b}-i{\bf c}$, where ${\bf b},{\bf c}\in\H$ are uniquely 
determined by a given ${\bf a}\in\M_2$. 

With this identification, a stem function
is a map $F:U\mapsto\M_2$, where $U\subset\C$ is conjugate
symmetric, with the property  $F({\bar\zeta})=\overline{F(\zeta)}$ for all $\zeta\in U$. Explicitly, if 
$$
F(\zeta)=\left(\begin{array}{cc} f_1(\zeta) & f_2(\zeta) 
\\ f_3(\zeta) & 
f_4(\zeta)  \end{array}\right),\,\,
F(\zeta)^\sim=\left(\begin{array}{cc} \overline{f_4(\zeta)} & -\overline{f_3(\zeta)} 
\\ -\overline{f_2(\zeta)} & 
\overline{f_1(\zeta)}  \end{array}\right),\,\,
\zeta\in U,
$$ 
% misprint: "\overline{f_4(\zeta)}" is replaced by  %"\overline{f_4(\zeta)}" 
then $F(\zeta)=F_1(\zeta)+iF_2(\zeta)$, with $F_1(\zeta)=F_1(\bar{\zeta})$ and  $F_2(\zeta)=-F_2(\bar{\zeta})$, where
$$
F_1(\zeta)=\frac{1}{2}\left(\begin{array}{cc} f_1(\zeta)+
\overline{f_4(\zeta)} & f_2(\zeta)-\overline{f_3(\zeta)} 
\\ f_3(\zeta)-\overline{f_2(\zeta)} & 
f_4(\zeta)+\overline{f_1(\zeta)} \end{array}\right),\,\,\zeta\in U,
$$ 
and 
$$
F_2(\zeta)=\frac{1}{2i}\left(\begin{array}{cc} f_1(\zeta)-
\overline{f_4(\zeta)} & f_2(\zeta)+\overline{f_3(\zeta)} 
\\ f_3(\zeta)+\overline{f_2(\zeta)} & 
f_4(\zeta)-\overline{f_1(\zeta)} \end{array}\right),\,\,\zeta\in U.
$$ 
In addition,
$$
f_1(\zeta)+
\overline{f_4(\zeta)}=f_1(\bar{\zeta})+
\overline{f_4(\bar{\zeta})},\,\,f_2(\zeta)-\overline{f_3(\zeta)}=
f_2(\bar{\zeta)}-\overline{f_3(\bar{\zeta})},
$$
$$
f_1(\zeta)-\overline{f_4(\zeta)}=-f_1(\bar{\zeta})+\overline{f_4(\bar{\zeta})},\,f_2(\zeta)+\overline{f_3(\zeta)}=
-f_2(\bar{\zeta})-\overline{f_3(\bar{\zeta})}.
$$
Consequently, $f_1(\zeta)=\overline{f_4(\bar{\zeta})}$, and
$f_2(\zeta)=-\overline{f_3(\bar{\zeta})},$ and so
$$
F(\zeta)=\left(\begin{array}{cc} f_1(\zeta) & f_2(\zeta) 
\\ -\overline{f_2(\bar{\zeta})} & 
\overline{f_1(\bar{\zeta})}  \end{array}\right),\,\,\zeta\in U,
$$
showing that every stem function is skew symmetric conjugate, via (\ref{func_scs}). 

Conversely, when $F$ is given by (\ref{func_scs}), so it is 
skew symmetric conjugate, we have
$$
F^\sim(\zeta)=\left(\begin{array}{cc} f_1(\bar{\zeta}) & f_2(\bar{\zeta)} 
\\ -\overline{f_2(\zeta)} & 
\overline{f_1(\zeta)}  \end{array}\right),\,\,\zeta\in U.
$$
Setting $F_1(\zeta)=(1/2)(F(\zeta)+F^\sim(\zeta))$ and 
$F_2(\zeta)=(1/2i)(F(\zeta)-F^\sim(\zeta))$, which are clearly
$\H$-valued functions, and $F(\zeta)=F_1(\zeta)+iF_2(\zeta)$, a direct computation shows that
$F_(\bar{\zeta})=F_1(\zeta)-iF_2(\zeta)$, showing that $F$ is a 
stem function. Some easy details are left to the reader.

As the term ''stem function`` is currently used in literature,
from now on we shall designate
a skew symmetric function as a stem function. Nevertheless, we shall use the definition of the skew symmetric function rather
than that equivalent of stem function, which is more appropriate to our framework. 

Finally, note that a stem function is not necessarily $\H$-valued. Using the notation from above, the stem function $F$ is
$\H$-valued if and only if $f_1(\bar{\zeta})=f_1(\zeta)$ and $f_2(\bar{\zeta})=-f_2(\zeta)$ for all $\zeta\in U$. 
\end{Rem}

% Rem 6

\begin{Rem}\label{calc_func}\rm With the notation from 
Definition \ref{canon_ev}, and because for each $\z\in\C^2$ the operator $Q(\z)$ is normal on the Hilbert space $\C^2$, we have a 
direct sum decomposition $\C^2=N_+(\z)\oplus N_-(\z)$, where 
$N_\pm(\z)=\{{\bf w}\in\C^2; Q(\z){\bf w}=s_\pm(\z){\bf w}\}.$
The projections $E_\pm(\z)$ of $\C^2$ onto $N_\pm(\z)$ are given by
$E_\pm(\z){\bf w}=\langle {\bf w},\nu_\pm(\z)\rangle\nu_\pm(\z),\,
{\bf w}\in\C^2$.

For every function $f:\sigma(Q({\bf z}))\mapsto\C$ we may define the operator
\begin{equation}\label{gen_func_calc}
f(Q({\bf z})){\bf w}=f(s_{+}({\bf z}))\langle{\bf w},\nu_{+}(\z)\rangle \nu_{+}(\z)+f(s_{-}({\bf z}))\langle{\bf w},\nu_{-}(\z)\rangle \nu_{-}(\z),
\end{equation}
where ${\bf w}\in\C^2$ is arbitrary, with a slight but traditional abuse of notation. 

We note  that 
formula (\ref{gen_func_calc}) is a particular case of the functional calculus given by the spectral
theorem for compact normal operators (see for instance \cite{GoGo}, Section VIII.2, or \cite{Kub}, Chapter 3). 

In particular, for $f(\zeta)=(\lambda-\zeta)^{-1},\,\zeta\in\sigma(Q(\z))),\,\lambda\neq\zeta$, we have
$$
(\lambda{\bf I}-Q(\z))^{-1}=(\lambda-s_{+}({\bf z}))^{-1}E_+(\z)+
(\lambda-s_{-}({\bf z}))^{-1}E_-(\z),\,\,\lambda\notin\sigma(Q(\z)),
$$
a formula to be later used.

More generally, 
for a fixed $\z=(z_1,z_2)\in\C^2$ and a function $F:\sigma(Q({\bf z}))\mapsto
\M_2$, we may define the operator (in fact a matrix with respect
to the canonical basis of $\C^2$) by the formula
\begin{equation}\label{gen_func_calc_vect}
F(Q({\bf z})){\bf w}=F(s_{+}({\bf z}))\langle{\bf w},\nu_{+}(\z)\rangle \nu_{+}(\z)+F(s_{-}({\bf z}))\langle{\bf w},\nu_{-}(\z)\rangle \nu_{-}(\z),
\end{equation}
where ${\bf w}\in\C^2$ is arbitrary, as an extension of 
(\ref{gen_func_calc}). Using the notation from Definition
\ref{canon_ev}, formula (\ref{gen_func_calc_vect}) can be also written as 
$$
F(q){\bf w}=F(s_{+}(q))\langle{\bf w},\nu_{+}(q)\rangle \nu_{+}(q)+F(s_{-}(q))\langle{\bf w},\nu_{-}(q)\rangle \nu_{-}(q),\,\,
{\bf w}\in\C^2,
$$
for each $F:\sigma(q)\mapsto\M_2$. In addition, if $q=Q(\z)$, one can use the notation $N_\pm(q)=N_\pm(\z)$ and 
$E_\pm(q)=E_\pm(\z)$. In fact, if  $U\subset\C$ is conjugate symmetric, the formula from above leads to a function $F:U_\H\mapsto\M_2$ (keeping, as usual, the same notation). Finally, 
when $q=s{\bf I},\,s\in\R$, then $F(q)=F(s){\bf I}$, via Lemma
\ref{roots}(c).  
\end{Rem}

% Thm 1

\begin{Thm}\label{sym_spec} Let $U\subset\C$ be a conjugate symmetric subset, and let $F:U\mapsto\M_2$.
The matrix $F(q)$ is a quaternion for all $q\in U_\H$ if and only if  $F$ is a stem function.
\end{Thm}

{\it Proof}\, We fix a point $\zeta\in U$. As $\bar{\zeta}\in U$, we may assume, with no loss of generality, that $\Im\zeta\ge0$. 

{\bf Case 1} We assume that $\Im\zeta>0$, and choose a quaternion $q\in U_\H$ with $\sigma(q)=\{\zeta,\bar{\zeta}\}$. 
Writing  $q=Q(\z)$ with ${\bf z}=(z_1,z_2)\in\C^2$,  because 
$\Im\zeta>0$, we may assume  $z_2\neq0$,
via Remark \ref{samesp}. Let $\nu_\pm(\z)$ be the canonical eigenvectors of $Q(\z)$, given by  Definition \ref{canon_ev}.
We also have $s_+(\z)=\zeta,\,s_-(\z)=\bar{\zeta}$. 

We show first that $F(Q({\bf z}))\in\H$ if and only if 
\begin{equation}\label{gen_H_cond}
F(s_{+}({\bf z}))\nu_{+}(\z)=F^\sim(s_{-}({\bf z}))\nu_{+}(\z).
\end{equation}

Let us look for the matrix form of $F(Q({\bf z}))$. To simplify
the computation, we set $s_\pm=s_\pm(\z), 
F_\pm=F(q_{\pm}({\bf z})), \nu_{\pm}(\z)=(\nu_{{\pm}1},\nu_{{\pm}2})$, and fix a ${\bf w}=(w_1,w_2)\in\C^2$. 

Note that
$$
(w_1\overline{\nu_{+1}}+w_2\overline{\nu_{+2}})\left( \begin{array}{cc}  \nu_{+1} \\ \nu_{+2} \end{array}\right)=
\left( \begin{array}{cc}  \vert\nu_{+1}\vert^2 & 
\nu_{+1}\overline{\nu_{+2}} \\ \overline{\nu_{+1}} \nu_{+2} & \vert\nu_{+2}\vert^2\end{array}\right)\left( \begin{array}{cc} w_1 \\ w_2 \end{array}\right),
$$
\begin{equation}\label{reverse}
{ }
\end{equation}
$$
(w_1\overline{\nu_{-1}}+w_2\overline{\nu_{-2}})\left( \begin{array}{cc}  \nu_{-1} \\ \nu_{-2} \end{array}\right)=
\left( \begin{array}{cc}  \vert\nu_{+2}\vert^2 & 
-\nu_{+1}\overline{\nu_{+2}}\\ -\overline{\nu_{+1}} \nu_{+2} & 
\vert\nu_{+1}\vert^2\end{array}\right) \left( \begin{array}{cc} w_1 \\ w_2 \end{array}\right),
$$
via Lemma \ref{nu_proper}.

Setting
$$
F_+=\left( \begin{array}{cc} f^+_{11} & f^+_{12} \\ f^+_{21} &  f^+_{22}\end{array}\right),\,\,
F_-=\left( \begin{array}{cc} f^-_{11} & f^-_{12} \\ f^-_{21} &  f^-_{22}\end{array}\right), 
$$
formulas (\ref{gen_func_calc_vect}) and (\ref{reverse}) lead to

$$
F(Q({\bf z}))=\left( \begin{array}{cc} f^+_{11} & f^+_{12} \\ f^+_{21} &  f^+_{22}\end{array}\right)
\left( \begin{array}{cc}  \vert \nu_{+1}\vert^2 & 
\nu_{+1}\overline{\nu_{+2}} \\ \overline{\nu_{+1}} \nu_{+2} & \vert\nu_{+2}\vert ^2\end{array}\right)+
$$
$$
\left( \begin{array}{cc} f^-_{11} & f^-_{12} \\ f^-_{21} &  f^-_{22}\end{array}\right)\left( \begin{array}{cc}  \vert\nu_{+2}\vert^2 & 
-\nu_{+1}\overline{\nu_{+2}}\\ -\overline{\nu_{+1}} \nu_{+2} & 
\vert\nu_{+1}\vert^2\end{array}\right).
$$

Note that
$$
\left( \begin{array}{cc}  \vert \nu_{+1}\vert^2 & 
\nu_{+1}\overline{\nu_{+2}} \\ \overline{\nu_{+1}} \nu_{+2} & \vert\nu_{+2}\vert ^2\end{array}\right)^{\sim}=\left( \begin{array}{cc}  \vert\nu_{+2}\vert^2 & 
-\nu_{+1}\overline{\nu_{+2}}\\ -\overline{\nu_{+1}} \nu_{+2} & 
\vert\nu_{+1}\vert^2\end{array}\right).
$$

Therefore, the matrix $ F(Q(\bf{z}))$ is a quaternion if and only if 
$$
\left(\left( \begin{array}{cc} f^+_{11} & f^+_{12} \\ f^+_{21} &  f^+_{22}\end{array}\right)-\left( \begin{array}{cc} f^-_{11} & f^-_{12} \\ f^-_{21} &  f^-_{22}\end{array}\right)^{\sim}\right)
\left( \begin{array}{cc}  \vert \nu_{+1}\vert^2 & 
\nu_{+1}\overline{\nu_{+2}} \\ \overline{\nu_{+1}} \nu_{+2} & \vert\nu_{+2}\vert ^2\end{array}\right)
$$
is a quaternion, via the properties of the skew complex conjugation.  Equivalently, $F(Q({\bf z}))$ is a quaternion if and only if the matrix
$$
\left(\begin{array}{cc} (f^+_{11}-\overline{f^-_{22}})\vert\nu_{+1}\vert^2+(f^+_{12}+\overline{f^-_{21}})\overline{\nu_{+1}}\nu_{+2} & (f^+_{11}-\overline{f^-_{22}})\nu_{+1}\overline{\nu_{+2}}+(f^+_{12}+\overline{f^-_{21}})\vert\nu_{+2}\vert^2 \\ (f^+_{21}+\overline{f^-_{12}})\vert\nu_{+1}\vert^2+(f^+_{22}-\overline{f^-_{11}})\overline{\nu_{+1}}\nu_{+2} & (f^+_{21}+\overline{f^-_{12}})\overline{\nu_{+1}}\nu_{+2}+(f^+_{22}-\overline{f^-_{11}})\vert\nu_{+2}\vert^2 
\end{array}\right)$$
is a quaternion, that is, 
$$
(f^+_{11}-\overline{f^-_{22}})\vert\nu_{+1}\vert^2+
(f^-_{11}-\overline{f^+_{22}})\vert\nu_{+2}\vert^2+
(f^+_{12}-f^-_{12}-\overline{f^+_{21}}+\overline{f^-_{21}})\overline{\nu_{+1}}\nu_{+2}=0
$$
\begin{equation}\label{sq-cond2}
{}
\end{equation}
$$
(f^-_{12}+\overline{f^+_{21}})\vert\nu_{+1}\vert^2+
(f^+_{12}+\overline{f^-_{21}})\vert\nu_{+2}\vert^2+
(f^+_{11}-f^-_{11}+\overline{f^+_{22}}-\overline{f^-_{22}})
\nu_{+1}\overline{\nu_{+2}}=0.
$$
\medskip

Let us make condition (\ref{sq-cond2}) more explicit.
To have a notation even simpler, we set $x:=f^+_{11}-
\overline{f^-_{22}},\,y:=f^-_{11}-\overline{f^+_{22}},\,
u:=f^+_{12}+\overline{f^-_{21}}, 
v=f^-_{12}+\overline{f^+_{21}}.$ We also put 
$a:=\nu_{+1},\,b:=\nu_{+2},\,c:=a\bar{b}$. Equations (\ref{sq-cond2}) become
$$
\vert a\vert^2x+\vert b\vert^2y+\bar{c}(u-v)=0,\,\,
\vert a\vert^2v+\vert b\vert^2u+c(x-y)=0.
$$
Equivalently,
$$
\left(\begin{array}{cc} \vert a\vert^2 & \vert b\vert^2 \\
c & -c\end{array}\right)\left(\begin{array}{cc} x
\\ y\end{array}\right)=\left(\begin{array}{cc} -\bar{c} & \bar{c} \\-\vert b\vert^2 & -\vert a\vert^2\end{array}\right)
\left(\begin{array}{cc} u \\ v \end{array}\right).
$$
As we have 
$$
\left(\begin{array}{cc} \vert a\vert^2 & \vert b\vert^2 \\
c & -c\end{array}\right)^{-1}=\frac{1}{c(\vert a\vert^2 + \vert b\vert^2)}\left(\begin{array}{cc} c & \vert b\vert^2\\ c &  -\vert a\vert^2\end{array}\right),
$$
we deduce that
$$
\left(\begin{array}{cc} x
\\ y\end{array}\right)=\frac{1}{c(\vert a\vert^2 + \vert b\vert^2)}\left(\begin{array}{cc} -\vert a\vert^2\vert b\vert^2
-\vert b\vert^4 & 0 \\ 0 & \vert a\vert^2\vert b\vert^2
+\vert a\vert^4\end{array}\right)\left(\begin{array}{cc} u \\ v \end{array}\right), 
$$
whence
$$
x=-\frac{b}{a}u,\,\,y=\frac{\bar{a}}{\bar{b}}v.
$$
Because 
$$
a=\frac{z_2}{\sqrt{\vert z_2\vert^2+\vert s_{+}(\z)-z_1\vert^2}},\,\,
b=\frac{s_{+}(\z)-z_1}{\sqrt{\vert z_2\vert^2+\vert s_{+}(\z)-z_1\vert^2}},
$$
we obtain
$$
f^+_{11}-\overline{f^-_{22}}=\frac{z_1-s_{+}(\z)}{z_2}
(f^+_{12}+\overline{f^-_{21}}),
$$
$$
f^-_{11}-\overline{f^+_{22}}=\frac{\bar{z}_2}{s_{-}(\z)-
\bar{z}_1}(f^-_{12}+\overline{f^+_{21}}).
$$

Therefore,
$$
f^+_{11}z_2+f^+_{12}(s_{+}(\z)-z_1)=\overline{f^-_{22}}z_2-\overline{f^-_{21}}(s_{+}(\z)-z_1),
$$
$$
f^+_{21}z_2+f^+_{22}(s_{+}(\z)-z_1)=-\overline{f^+_{12}}z_2+\overline{f^-_{11}}(s_{+}(\z)-z_1),
$$
and so
$$
\left(\begin{array}{cc} f^+_{11} & f^+_{12}\\
f^+_{21} & f^+_{22}\end{array}\right)
\left(\begin{array}{cc} z_2 \\ s_{+}(\z)-z_1 \end{array}\right)=
\left(\begin{array}{cc} \overline{f^-_{22}} & -\overline{f^-_{21}}\\ -\overline{f^- _{12}} & \overline{f^-_{11}}\end{array}\right)\left(\begin{array}{cc} z_2 \\ s_{+}(\z)-z_1 \end{array}\right),
$$
which is (\ref{gen_H_cond}),  modulo a multiplicative constant.
In other words, $F(Q(\z))\in\H$ if and only if (\ref{gen_H_cond})
holds. 

Next we apply Remark \ref{samesp}. Because $\Im\zeta>0$, we define
$$
\z_u=(\Re\zeta+i\sqrt{(\Im\zeta)^2-\vert u\vert^2},u)\in\C^2,
$$
for some $u\in\C$ with $0\neq\vert u\vert<\Im\zeta$. Let also
$\nu_+(\z_u)$ be the corresponding canonical eigenvector of $Q(\z_u)$. Because
$s_+(\z_u)=\zeta=s_+(\z)$, the previous argument shows that
$F(Q(\z_u))$ is a quaternion if and only if $F(\zeta)\nu_+(\z_u)=
F^\sim(\bar{\zeta})\nu_+(\z_u)$. Now, there are pairs 
$u_1,u_2$ such that $\nu_+(\z_{u_1}),\,\nu_+(\z_{u_2})$ are 
linearly independent. For instance, if $0< u_1=-u_2<\Im(\zeta)$, the vectors 
$$
(u_1,i(\Im\zeta-\sqrt{(\Im\zeta)^2-u_1^2})),\,(u_2,i(\Im\zeta-\sqrt{(\Im\zeta)^2-u_2^2}))
$$
are linearly independent. Indeed, if not, we 
would have 
$$
\Im\zeta-\sqrt{(\Im\zeta)^2-u_1^2}=-\Im\zeta+
\sqrt{(\Im\zeta)^2-u_1^2}),
$$
implying $u_1=0$, which is not possible. Having two  linearly independent vectors $\nu_+(\z_{u_1}),\,\nu_+(\z_{u_2})$, with the hypothesis $F(Q(\z_{u_1})), F(Q(\z_{u_2}))\in\H$ we deduce  
that the equality $F(\zeta)=F^\sim(\bar{\zeta})$ holds. 

If $\Im\zeta<0$, a similar argument shows that 
$F(\zeta)=F^\sim(\bar{\zeta})$ (see also Remark \ref{sq_gen}). In particular, if $F(Q(\z))\in\H$
for $Q(\z)\in U_\H$, the equality $F(\zeta)=F^\sim(\bar{\zeta})$ is true for all $\zeta\in U$ with $\Im\zeta\neq0$.

{\bf Case 2} We assume that $\zeta\in U\cap\R$, and put 
$x=\zeta$. When $\z=(x,0)$ with $x\in\R$, we have 
$Q(\z))=x\bf I$ and so $F(Q(\z))=F(x)\bf I$. In this case, it is
obvious that 
$F(x)\bf I$ is a quaternion if and only if  $F(x)=F^\sim(x)$. 
 
Consequently, if $F(q)\in\H$ for all $q\in U_\H$,
the function $F$ is a stem one.

{\bf Final Case} To finish the proof, we have to show that if $F$ is a stem function, we must have $F(q)\in\H$ for all $q\in U_\H$.

Fixing a point $q=Q(\z)$ with $\z=(z_1,z_2)$, assuming $z_2\neq0$,
if $\zeta=s_+(\z)$, condition (\ref{gen_H_cond}) cleraly holds,
so $F(q)\in\H$. 

Next assume that $z_2=0$ but $\Im z_1\neq0$. In this case, using
the notation from Definition \ref{canon_ev}, we 
have $s_+(q)=z_1,s_-(q)=\bar{z}_1, \nu_+(q)=(1,0), \nu_-(q)=(0,1)$, and with $F$ as in formula (\ref{func_scs}), formula 
(\ref{gen_func_calc_vect}) leads to
$$
F(q){\bf w}=
F(z_1)\langle{\bf w},\nu_{+}(q)\rangle \nu_{+}(q)+
F^\sim(\bar{z}_1)\langle{\bf w},\nu_{-}(q)\rangle \nu_{-}(q)=
$$
$$
\left(\begin{array}{cc} f_1(\bar{z}_1) & f_2(z_1) 
\\ -\overline{f_2(\bar{z}_1)} & 
\overline{f_1(\bar{z}_1)}  \end{array}\right)\left(\begin{array}{cc} w_1 \\ 0\end{array}\right)+
\left(\begin{array}{cc} f_1(\bar{z}_1) & f_2(\bar{z}_1) 
\\ -\overline{f_2(z_1)} & 
\overline{f_1(z_1)}\end{array}\right)\left(\begin{array}{cc} 0 
\\ w_2\end{array}\right)= 
$$ 
$$
\left(\begin{array}{cc} f_1(z_1) & f_2(\bar{z}_1) 
\\ -\overline{f_2(\bar{z}_1)} & 
\overline{f_1(z_1)}  \end{array}\right)\left(\begin{array}{cc} w_1 \\ w_2\end{array}\right)
$$ 
for all ${\bf w}=(w_1,w_2)\in\C^2$, which shows that $F(q)$ is a
quaternion.   
 
Finally, if  $z_2=0$ and $\Im z_1=0$ the assertion is given by
Case 2, from above.

% Rem 7

\begin{Rem}\label{sq_gen}\rm 
 There is a parallel treatment 
 of Case 1 from the previous proof, leading to a  similar statement.  We first introduce the function
$$
G:=\left(\begin{array}{cc} f_{21} & f_{22} \\
-f_{11} & -f_{12}\end{array}\right):U\mapsto\M_2, 
$$
whose entries are given by the original function
$$
F=
\left( \begin{array}{cc} f_{11} & f_{12} \\ f_{21} & f_{22} \end{array}\right):U\mapsto\M_2.
$$
Then we prove that
\begin{equation}\label{gen_H_cond3}
G(s_{-}({\bf z}))\nu_{-}(\z)=G^\sim(s_{+}({\bf z}))\nu_{-}(\z).
\end{equation}
To infer (\ref{gen_H_cond3}), we use conditions (\ref{sq-cond2}),
written as 
$$
(f^+_{11}-\overline{f^-_{22}})\vert\nu_{-2}\vert^2+
(f^-_{11}-\overline{f^+_{22}})\vert\nu_{-1}\vert^2-
(f^+_{12}-f^-_{12}-\overline{f^+_{21}}+\overline{f^-_{21}})\overline{\nu_{-1}}\nu_{-2}=0.
$$
\begin{equation}\label{sq-cond3}
{}
\end{equation}
$$
(f^-_{12}+\overline{f^+_{21}})\vert\nu_{-2}\vert^2+
(f^+_{12}+\overline{f^-_{21}})\vert\nu_{-1}\vert^2-
(f^+_{11}-f^-_{11}+\overline{f^+_{22}}-\overline{f^-_{22}})
\nu_{-1}\overline{\nu_{-2}}=0.
$$
As in the proof of Case 1, we obtain the equality 
\begin{equation}\label{sq-cond4}
\left(\begin{array}{cc} -\overline{f^+_{12}} & 
\overline{f^+_{11}}\\
-\overline{f^+_{22}} & \overline{f^+_{21}}\end{array}\right)
\left(\begin{array}{cc} z_2 \\ s_{-}(\z)-z_1 \end{array}\right)=
\left(\begin{array}{cc} f^-_{21} & f^-_{22}\\ -f^-_{11} & 
-f^-_{12}\end{array}\right)\left(\begin{array}{cc} z_2 \\ s_{-}(\z)-z_1 \end{array}\right),
\end{equation}
which is relation (\ref{gen_H_cond3}), modulo a multiplicative
factor. Using, as above, Remark \ref{samesp}, we deduce the equality $G(s_-(\z))=G^\sim(s_+(\z))$, which is equivalent to  $F(s_-(\z))=F^\sim(s_+(\z))$.

We have, in fact, an equivalent formulation of Theorem
\ref{sym_spec}, asserting that $F(q)\in\H$ for all $q\in U_\H$ if and only if $G(\zeta)=G^\sim(\bar{\zeta})$ for all $\zeta\in U$, 
that is, if and only if $G$ is a stem function. 
\end{Rem}

% Cor 1

\begin{Cor}\label{sym_spec1}  Let $U\subset\C$ be a conjugate symmetric subset, and let $f:U\mapsto\C$.
We have $f(q)\in\H$ for all $q\in U_\H$ if and only if 
$f(\zeta)=\overline{f(\bar{\zeta})}$ for all $\zeta\in U$.
\end{Cor}

{\it Proof.}\, We apply Theorem \ref{sym_spec} to the function
$F=f{\bf I}:U\mapsto\M_2$. This function is  a stem one if and only if $f(\zeta)=\overline{f(\bar{\zeta})}$
for all $\zeta\in U$.

% Cor 2

\begin{Cor}\label{sym_spec2}  Let $U\subset\C$ be an open conjugate symmetric subset, and let $F:U\mapsto\H$. Then we
have $F(q)\in\H$ for all $q\in U_\H$ if and only if $F(\zeta)=F(\bar{\zeta})$ for all $\zeta\in U$.
\end{Cor}

{\it Proof.}\, The property $F:U\mapsto\H$ implies that
$F^\sim=F$.  Therefore, $F$ is a stem function
if and only if $F(\zeta)=F(\bar{\zeta})$ for all $\zeta\in U$.

% Rem 7

\begin{Rem}\label{extension}\rm Let $U\subset\C$ be a conjugate
symmetric set, and let $F:U\mapsto\M_2$ be a  stem function. The formula
$$
F(q){\bf w}=F(s_{+}(q))\langle{\bf w},\nu_{+}(q)\rangle \nu_{+}(q)+F(s_{-}(q))\langle{\bf w},\nu_{-}(q)\rangle \nu_{-}(q),
$$
where $q\in U_\H$ and ${\bf w}\in\C^2$ are arbitrary, is an 
''extension`` of the function $F$ to $U_\H$, in a sense to
be specified.  Note that  we have 
an embedding $U\ni\zeta\mapsto q_\zeta:=Q((\zeta,0))\in \H$,
which is the restriction of an $\R$-linear map isometry. In fact, writing $\zeta=x+iy$, with $x,y\in\R$ unique, we have 
$$
q_\zeta=\left(\begin{array}{cc} x+iy & 0 \\ 0 & x-iy\end{array}\right)=
x{\bf I}+y{\bf J},
$$ 
allowing us to identify 
the set $U$ with the set 
$$U_{\bf J}:=\{q_\zeta;\zeta\in U\}=\{x{\bf I}+y{\bf J}; x+iy\in U\}\subset\H.
$$ 
Because $F$ is a stem function, we must have
$$
F(\zeta)=\left(\begin{array}{cc} f_1(\zeta) & f_2(\zeta) 
\\ -\overline{f_2(\bar{\zeta})} & 
\overline{f_1(\bar{\zeta})}  \end{array}\right),\,\,\zeta\in U,
$$
given by formula (\ref{func_scs}).

As we have $\sigma(q_\zeta)=\{\zeta,\bar{\zeta}\}$, it follows from Lemma \ref{roots}
that $s_+(q_\zeta)=\zeta,\,s_-(q_\zeta)=\bar{\zeta},\,\nu_+(q_\zeta)=(1,0),\,\nu_-(q_\zeta)=(0,1)$, and hence, as in the proof of Theorem \ref{sym_spec}, 
$$
F(q_\zeta)\left(\begin{array}{cc} w_1 \\ w_2 \end{array}\right)=F(\zeta)\left(\begin{array}{cc} w_1 \\ 0 \end{array}\right)+F(\bar{\zeta})\left(\begin{array}{cc} 0 \\ w_2 \end{array}\right)=
$$
$$
\left(\begin{array}{cc} f_1(\zeta) & f_2(\bar{\zeta}) 
\\ -\overline{f_2(\bar{\zeta})} & \overline{f_1(\zeta)}\end{array}\right)
\left(\begin{array}{cc} w_1 \\ w_2 \end{array}\right)
$$
for all $(w_1,w_2)\in\C^2$,  so
$$
F(q_\zeta)=\left(\begin{array}{cc} f_1(\zeta) & f_2(\bar{\zeta}) 
\\ -\overline{f_2(\bar{\zeta})} & \overline{f_1(\zeta)}\end{array}\right),\,\,\zeta\in U.
$$

Let $\mathcal{F}_s(U,\M_2)=\{F:U\mapsto\M_2;
F(\bar{\zeta})=F^\sim(\zeta),\,\zeta\in U\}$, which is an
$\R$-algebra of $\M_2$-valued functions on $U$, consisting of all
stem functions on $U$. Let also 
$\mathcal{F}(U,\H)=\{G:U\mapsto\H\}$, which an $\R$-algebra of $\H$-valued functions on $U$. Setting $\kappa(\zeta)=q_\zeta,\,\zeta\in U$,
we have an injective unital morphism of $\R$-algebras given by
$\mathcal{F}_s(U,\M_2)\ni F\mapsto F\circ\kappa\in \mathcal{F}(U,\H)$. 
Therefore, the map $U_\H\ni q \mapsto F(q)\in\H$ given by
(\ref{gen_func_calc_vect}), which extends the map
$q_\zeta\mapsto F(q_\zeta)$, may be also regarded as an ''extension`` of $F\in \mathcal{F}_s(U,\M_2)$ (modulo the map
$\kappa$). Note also that
the function $U_\H\ni q \mapsto F(q)\in\H$ is uniquely determined by the function $U\ni\zeta\mapsto F(\zeta)\in\M_2$, when the latter is a stem function.

In particular, if $f\in\mathcal{F}_s(U):
=\{g:U\mapsto\C;g(\bar{\zeta})=
\overline{g(\zeta)},\,\zeta\in U\}$, then 
$$
f(Q(\zeta,0))=\left(\begin{array}{cc} f(\zeta) & 0 \\ 0 &
\overline{f(\zeta)}\end{array}\right)=Q((f(\zeta),0)),\,\,\zeta\in U.
$$
 
\end{Rem}

\section{Analytic Functional Calculus for Quaternions}
\label{AFCQ}

Regarding, as before, the quaternions as normal operators,
we now investigate some consequences of their analytic functional calculus, in the classical sense (see \cite{DuSc}, Section VII.3, for details). The frequent use of 
various versions of the Cauchy formula is simplified by adopting the following definition. Let $U\subset\C$ be open. An open subset $\Delta\subset U$ will be called a {\it Cauchy domain} (in $U$) if $\Delta\subset\bar{\Delta}\subset U$ and the boundary $\partial\Delta$ of $\Delta$ consists of a finite family of closed curves, piecewise smooth, positively oriented.
Note that a Cauchy domain is bounded but not necessarily connected.

% Lem 3

\begin{Lem}\label{vect_fc} Let $U\subset\C$ be a conjugate symmetric open set 
and let $F:U\mapsto\M_2$ be an analytic function. For every $q\in U_\H$ we set 
\begin{equation}\label{Cauchy_vect}
F_\H(q)=\frac{1}{2\pi i}\int_\Gamma F(\zeta)(\zeta{\bf I}-q)^{-1}d\zeta,
\end{equation}
where $\Gamma$ is the boundary of a Cauchy domain in $U$ containing the spectrum $\sigma(q)$. Then we have  $F_\H(q)\in\H$  for all  $q\in U_\H$ if and only if $F$ is a stem function. 
\end{Lem} 

{\it Proof.\,} We first assume that $q\notin\R{\bf I}$.  If 
$\sigma(q)=\{s_+,s_-\}$ with $s_\pm=s_\pm(q)$, the points $s_+,s_-$ are distinct and  not real, by Lemma \ref{roots}. We fix an $r>0$ sufficiently small such that,
setting $D_\pm:=\{\zeta\in U;\vert\zeta-s_\pm\vert\le r\}$, we
have 
$D_\pm\subset U$ and $D_+\cap D_-=\emptyset$. Then 
$$
F_\H(q)=\frac{1}{2\pi i}\int_{\Gamma_+} F(\zeta)(\zeta{\bf I}-q)^{-1}d\zeta+\frac{1}{2\pi i}\int_{\Gamma_-} F(\zeta)(\zeta{\bf I}-q)^{-1}d\zeta,
$$
where $\Gamma_\pm$ is the boundary of $D_\pm$. We may write 
$F(\zeta)=\sum_{k\ge0}(\zeta-s_+)^kA_k$ with $\zeta\in D_+$,
$A_k\in\M_2$ for all $k\ge0$, as a uniformly convergent series.
Similarly, $F(\zeta)=\sum_{k\ge0}(\zeta-s_-)^kB_k$ with $\zeta\in D_-$, $B_k\in\M_2$ for all $k\ge0$, as a uniformly convergent series. 
 
Note that
$$
\frac{1}{2\pi i}\int_{\Gamma_+} F(\zeta)(\zeta{\bf I}-q)^{-1}d\zeta=\sum_{k\ge0}\left(A_k\frac{1}{2\pi i}\int_{\Gamma_+} (\zeta-s_+)^k(\zeta{\bf I}-q)^{-1}d\zeta\right)=A_0E_+,
$$
because the integral
$$
\frac{1}{2\pi i}\int_{\Gamma_+} (\zeta-s_+)^k(\zeta{\bf I}-q)^{-1}d\zeta
$$
is equal to $E_+:=\frac{1}{2\pi i}\int_{\Gamma_+} (\zeta{\bf I}-q)^{-1}d\zeta$ when $k=0$, which is the projection of $\C^2$ onto the space
$N_+:=\{{\bf v}; q{\bf v}=s_+{\bf v}\}$, and it is equal to 
$0$ when $k\ge 1$, because $(\zeta{\bf I}-q)^{-1}E_+=
(\zeta-s_+)^{-1}E_+$ (see Remark \ref{calc_func}).

Similarly
$$
\frac{1}{2\pi i}\int_{\Gamma_-} F(\zeta)(\zeta{\bf I}-q)^{-1}d\zeta=\sum_{k\ge0}\left(B_k\frac{1}{2\pi i}\int_{\Gamma_-} (\zeta-s_-)^k(\zeta{\bf I}-q)^{-1}d\zeta\right)=B_0E_-
$$
because, as above, the integral
$$
\frac{1}{2\pi i}\int_{\Gamma_-} (\zeta-s_-)^k
(\zeta{\bf I}-q)^{-1})d\zeta
$$
is equal to 
$E_-:=\frac{1}{2\pi i}\int_{\Gamma_-} (\zeta{\bf I}-q)^{-1}d\zeta$ when $k=0$, which is the projection of $\C^2$ onto the space
$N_-:=\{{\bf v}; q{\bf v}=s_-{\bf v}\}$, and it is equal to 
$0$ when $k\ge 1$, since $(\zeta{\bf I}-q)^{-1}E_-=
(\zeta-s_-)^{-1}E_-$. 
Note also that $E_+{\bf w}=\langle{\bf w},\nu_{+}(q)\rangle \nu_{+}(q)$, and  
$E_-{\bf w}=\langle{\bf w},\nu_{-}(q)\rangle \nu_{-}(q)$, 
all ${\bf w}\in\C^2$. Consequently,
$$
F_\H(q)=F(s_+)E_++F(s_-)E_-,
$$
and the right hand side of this equality coincides with formula
(\ref{gen_func_calc_vect}). 
 
Assume now that
$\sigma(q)=\{s\}$,  where  $s:=s_+=s_-\in\R$. We fix an $r>0$ such that  
$D:=\{\zeta\in U;\vert\zeta-s\vert\le r\}\subset U$. Then we have
$$
F_\H(q)=\frac{1}{2\pi i}\int_{\Gamma} F(\zeta)(\zeta{\bf I}-q)^{-1}d\zeta,
$$
where $\Gamma$ is the boundary of $D$. We write 
$F(\zeta)=\sum_{k\ge0}(\zeta-s)^kA_k$ with $\zeta\in D$,
$A_k\in\M_2$ for all $k\ge0$, as a uniformly convergent series.
Note also that
$$
\frac{1}{2\pi i}\int_{\Gamma} F(\zeta)(\zeta{\bf I}-q)^{-1}d\zeta=\sum_{k\ge0}\left(A_k\frac{1}{2\pi i}\int_{\Gamma} (\zeta-s)^k(\zeta{\bf I}-q)^{-1}d\zeta\right)=A_0,
$$
because the integral 
$$
\frac{1}{2\pi i}\int_{\Gamma} (\zeta-s)^k(\zeta{\bf I}-q)^{-1}d\zeta
$$
is equal to $\bf I$ when $k=0$, and equal to $0$ when
$k\ge 1$, since $(\zeta{\bf I}-q)^{-1}=(\zeta-s)^{-1} {\bf I}$. Consequently, $F_\H(q)=A_0=F(s){\bf I}$. 

In both situations, the matrix $F_\H(q)$ is equal to  the right hand side of formula (\ref{gen_func_calc_vect}). Therefore, we must have 
$F_\H(q)\in\H$ if and only if $F(s_+)=F^\sim(s_-)$, via Theorem 
\ref{sym_spec}. In other words, $F_\H(q)\in\H$ for all 
$q\in U_\H$ if and only if $F:U\mapsto \M_2$ is a stem function.  
 
% Rem 8

\begin{Rem}\label{bi_not}\rm It follows from the proof of the previous lemma
that the element $F_\H(q)$, given by formula (\ref{Cauchy_vect}), coincides with the element $F(q)$ given 
by (\ref{gen_func_calc_vect}). Nevertheless, we keep the 
notation $F_\H(q)$ whenever we want to emphasize that it is 
defined via (\ref{Cauchy_vect}).
\end{Rem}

% Cor 3

\begin{Cor}\label{vect_fc_cor} Let $U\subset\C$ be a conjugate symmetric open set 
and let $f:U\mapsto\C$ be an analytic function. For every $q\in U_\H$ we set 
\begin{equation}\label{Cauchy_vect1}
f_\H(q)=\frac{1}{2\pi i}\int_\Gamma f(\zeta)(\zeta{\bf I}-q)^{-1}d\zeta,
\end{equation}
where  $\Gamma$ is the boundary of a Cauchy domain in $U$ containing the spectrum $\sigma(q)$. Then we have  $f_\H(q)\in\H$ if and only if
$f(s_+(q))=\overline{f(s_-(q))}$ for all $q\in U_\H$. 
\end{Cor} 

{\it Proof.}\, The assertion is a direct consequence of 
Lemma \ref{vect_fc}, applied to the the function $f\bf{I}$. 

% Rem 10

\begin{Rem}\label{ext}\rm Let $U\subset\C$ be open and conjugate symmetric. As already seen, for every point $\zeta\in U$ the quaternion $Q((\zeta,0))$ is an element of $U_\H$ because its spectrum
equals the set $\{\zeta,\bar{\zeta}\}$. According to Corollary
\ref{vect_fc_cor}, 
an analytic function $f:U\mapsto\C$ has the property $f_\H(q)\in 
\H$ for all $q\in U_\H$ if and only if $f(\bar{\zeta})=\overline{f(\zeta)}$ for all $\zeta\in U$. This shows, in particular, that the set $\mathcal{O}_s(U)$, consisting of all
analytic functions $f:U\mapsto\C$ with the property $f(\bar{\zeta})=\overline{f(\zeta)}$ for all $\zeta\in U$,
is compatible with the analytic functional calculus of 
the quaternions. Clearly, $\mathcal{O}_s(U)$ is a unital $\R$-subalgebra of the $\C$-algebra 
$\mathcal{O}(U)$ of all analytic functions in $U$.  

More generally, if $F\in\mathcal{O}(U,\M_2)$, where  
$\mathcal{O}(U,\M_2)$ is the $\C$-algebra of all $\M_2$-valued analytic functions in $U$,
we have the property $F_\H(q)\in\H$ for all $q\in U_\H$ if and only if $F$ is a stem function, by Lemma \ref{vect_fc}.  In this case, we have
$$
F(\zeta)=\left(\begin{array}{cc} f_1(\zeta) & f_2(\zeta) 
\\ -\overline{f_2(\bar{\zeta})} & 
\overline{f_1(\bar{\zeta})}  \end{array}\right),\,\,\zeta\in U,
$$
by formula (\ref{func_scs}). 

Let us denote by $\mathcal{O}_s(U,\M_2)$ the set of all analytic stem functions $F\in \mathcal{O}(U,\M_2)$, so $\mathcal{O}_s(U,\M_2)\subset\mathcal{F}_s(U,\M_2)$, 
where the latter is introduced in Remark \ref{extension}.

Note that if  $f_1\in\mathcal{O}_s(U)$ and  $f_2\in i\mathcal{O}_s(U)$,
we have
$$
F(\zeta)=\left(\begin{array}{cc} f_1(\zeta) & f_2(\zeta) 
\\ f_2(\zeta) & 
f_1(\zeta)  \end{array}\right)\in\mathcal{O}_s(U,\M_2).
$$

The converse is not true, in general. For instance, if $f(\zeta)=\zeta+i$ on $U=\C$, we have $\overline{f(\bar{\zeta})}=\zeta-i$, so 
$f\notin\mathcal{O}_s(U)$ but  
$$
F(\zeta)=\left(\begin{array}{cc}\zeta+i & 0
\\ 0 & 
\zeta-i \end{array}\right)\in\mathcal{O}_s(U,\M_2).
$$ 
 
It is easily seen that $\mathcal{O}_s(U,\M_2)$ is a unital 
$\R$-subalgebra of the $\C$-algebra $\mathcal{O}(U,\M_2)$, and 
it is also an  $\mathcal{O}_s(U)$-module.  

Finally, if $\Delta\subset\C$ is an open disk centered at $0$, each function $F\in\mathcal{O}_s(\Delta,\M_2)$ can be represented 
as a convergent series $F(\zeta)=\sum_{k\ge 0}a_k\zeta^k,\,\zeta\in\Delta$, with $a_k\in\H$ for all $k\ge 0$. 
 \end{Rem}

\begin{Def}\label{reg_func}\rm Let $\Omega\subset\H$ be a spectrally saturated open set, and 
let $U=\mathfrak{S}(\Omega)\subset\C$ (which is also open by
Remark \ref{consym}(4)). Then we put
$$
\mathcal{R}(\Omega)=\{f_\H; f\in\mathcal{O}_s(U)\},
$$
and
$$
\mathcal{R}(\Omega,\H)=\{F_\H; F\in\mathcal{O}_s(U,\M_2)\}.
$$
\end{Def} 

In fact, these are $\R$-linear spaces, having some important 
properties:

% Thm 2

\begin{Thm}\label{H_afc} Let $\Omega\subset\H$ be a spectrally
saturated open set, and let $U=\mathfrak{S}(\Omega)\subset\C$.
The space $\mathcal{R}(\Omega)$ is a unital commutative $\R$-algebra, the space $\mathcal{R}(\Omega,\H)$ is a right 
$\mathcal{R}(\Omega)$-module, and the map
$$
{\mathcal O}_s(U,\M_2)\ni F\mapsto F_\H\in\mathcal{R}(\Omega,\H)
$$ 
is a right module isomorphism. Moreover, for every polynomial\,\,$P(\zeta)=\break\sum_{n=0}^m a_n\zeta^n,\,\zeta\in\C$, with $a_n\in\H$  for all $n=0,1,\ldots,m$, we have  $P_\H(q)=\sum_{n=0}^m a_n q^n\in\H$ for all $q\in\H$.   
\end{Thm}

{\it Proof.\,} The $\R$-linearity of the maps
$$
{\mathcal O}_s(U,\M_2)\ni F\mapsto F_\H\in\mathcal{R}(\Omega,\H),\,
{\mathcal O}_s(U)\ni f\mapsto f_\H\in\mathcal{R}(\Omega),
$$
is clear. As the second one is also multiplicative follows from
the multiplicativiry of the analytic functional calculus at any 
point $q\in\H$.

In fact, we have a more general property, specifically $(Ff)_\H(q)=F_\H(q)f_\H(q)$ for all   $F\in {\mathcal O}_s(U,\M_2),\, f\in{\mathcal O}_s(U)$, and $q\in\Omega$. This follows from the equalities,
$$
(Ff)_\H(q)=\frac{1}{2\pi i}\int_{\Gamma_0} F(\zeta)f(\zeta)(\zeta{\bf I}-q)^{-1}d\zeta=
$$
$$
\left(\frac{1}{2\pi i}\int_{\Gamma_0} F(\zeta)(\zeta{\bf I}-q)^{-1}d\zeta\right)
\left(\frac{1}{2\pi i}\int_{\Gamma} f(\eta)(\eta{\bf I}-q)^{-1}d\eta\right)=F_\H(q)f_\H(q),
$$
obtained as in the classical case (see \cite{DuSc}, Section VII.3), which holds because $f$ is $\C$-valued and commutes with the quaternions in 
$\M_2$. Here $\Gamma,\,\Gamma_0$ are the boundaries of two Cauchy
domains $\Delta,\,\Delta_0$ respectively, such that $\Delta\supset
\bar{\Delta}_0$, and $\Delta_0$ contains $\sigma(q)$. 

Note that, in particular,  for every polynomial $P(\zeta)=\sum_{n=0}^m a_n\zeta^n$ with $a_n\in\H$  for all $n=0,1,\ldots,m$, we have  $P_\H(q)=\sum_{n=0}^m a_n q^n\in\H$ for all $q\in\H$. 

Another stated property is the injectivity of the map 
$$
{\mathcal O}_s(U,\M_2)\ni F\mapsto F_\H\in
\mathcal{R}(\Omega,\H)
$$ 
Indeed,  if the function $F_\H$ is null,
the function $U\ni\zeta\mapsto 
F_\H(q_\zeta)\in\H$ is null too, so $F\in{\mathcal O}_s(U,\M_2)$ should be null as well (see Remark \ref{extension}). 

 % Cor 4

\begin{Cor}\label{H_afcc} Let $\Omega\subset\H$ be a spectrally
saturated open set, and let $U=\mathfrak{S}(\Omega)\subset
\C$. The map
$$
\mathcal{O}_s(U)\ni f\mapsto f_\H\in\mathcal{R}(\Omega)
$$ 
is a unital $\R$-algebra isomorphism. Moreover, 

$(a)$ for every polynomial $p(\zeta)=\sum_{n=0}^m a_n\zeta^n$ with $a_n$ real for all $n=0,1,\ldots,m$, we have  $p_\H(q)=\sum_{n= 0}^m a_n q^n\in\H$ for all $q\in\Omega$;

$(b)$ if $f\in\mathcal{O}_s(U)$ has no zero in $U$, we 
have $(f_\H(q))^{-1}=f^{-1}_\H(q)$ for all $q\in\Omega$.   
\end{Cor}

The assertions are direct consequences of the previous proof.

% Cor 5

\begin{Cor}\label{disc} Let $r>0$ and let $U\supset\{\zeta\in\C;\vert\zeta\vert\le r\}$ be a conjugate symmetric open set. Then for every
$F\in\mathcal{O}_s(U,\M_2)$ one has
$$
F_\H(q)=\sum_{n\ge0}\frac{F^{(n)}(0)}{n!} q^n,\,\,\Vert q\Vert<r,
$$
where the series is absolutely convergent. 
\end{Cor}

{\it Proof.}\, Indeed, as we have $(\zeta{\bf I}-q)^{-1}=
\sum_{n\ge0}\zeta^{-n-1}q^n$, where the series is uniformly 
convergent on $\{\zeta;\vert\zeta\vert=r\}$ for each fixed $q$ with $\Vert q\Vert<r$, we infer that
$$
F_\H(q)=\frac{1}{2\pi i}\int_{\vert\zeta\vert=r}F(\zeta)
(\zeta{\bf I}-q)^{-1}d\zeta=
$$
$$
\sum_{n\ge0}\left(\frac{1}{2\pi i}\int_{\vert\zeta\vert=r}
\frac{F(\zeta)}{\zeta^{n+1}}d\zeta\right)q^n=
\sum_{n\ge0}\frac{F^{(n)}(0)}{n!} q^n.
$$

% Rem 10

\begin{Rem}\label{n-deriv}\rm For every function $F\in\mathcal{O}_s(U,\M_2)$, the derivatives $F^{(n)}$ also belong to 
$\mathcal{O}_s(U,\M_2)$, where $U\subset\C$ is a conjugate symmetric open set. To check this assertion, we note first 
that if $f\in\mathcal{O}(U)$, setting $f^*(\zeta):=\overline{f(\bar{\zeta})}$ for all $\zeta\in U$, we can easily
see that $f^*\in\mathcal{O}(U)$ and $(f^*)'(\zeta)=\overline{f'(\bar{\zeta})}$, for all $\zeta\in U$. 

Because
$$
F(\zeta)=\left(\begin{array}{cc} f_1(\zeta) & f_2(\zeta) 
\\ -f_2^*(\zeta) & 
f_1^*(\zeta)  \end{array}\right),\,\,\zeta\in U,
$$
we have
$$
F(\zeta)'=\left(\begin{array}{cc} f_1'(\zeta) & f_2'(\zeta) 
\\ -(f_2^*)'(\zeta) & 
(f_1^*)'(\zeta)  \end{array}\right),\,\,\zeta\in U,
$$
showing that  $F'\in\mathcal{O}_s(U,\M_2)$. 

Now fixing $F\in\mathcal{O}_s(U,\M_2)$, we may define its 
{\it extended derivatives} with respect to the quaternionic variable via the formula

\begin{equation}\label{Cauchy_vect_deriv}
F^{(n)}_\H(q)=\frac{1}{2\pi i}\int_\Gamma F^{(n)}(\zeta)(\zeta{\bf I}-q)^{-1}d\zeta,
\end{equation}
for the boundary $\Gamma$ of a Cauchy domain $\Delta\subset U$, $n\ge0$ an arbitrary integer, and $\sigma(q)\subset\Delta$. 

In particular, if $F\in\mathcal{O}_s(\Delta,\M_2)$, with $\Delta$
a disk centered at zero, and so we have a representation  as a convergent series $F(\zeta)=\sum_{k\ge0}a_k\zeta^k$  with coefficients in $\H$, then 
(\ref{Cauchy_vect_deriv}) gives the equality
$F'_\H(q)=\sum_{k\ge1}ka_kq^{k-1}$, which looks like a (formal) derivative of the function $F_\H(q)=\sum_{k\ge0}a_kq^{k}$. 
\end{Rem}

% Rem 11

\begin{Rem}\label{clas_C_ineg}\rm Let $U\subset\C$ be an arbitrary open set and let 
$F\in\mathcal{O}(U,\M_2)$ be an arbitrary analytic function.
Let also $\Delta_0,\Delta$ be open discs in $U$ with boundaries
$\Gamma_0,\Gamma$ respectively, such that
$\Delta_{0}\subset\bar{\Delta}_{0}\subset\Delta\subset\bar{\Delta}\subset U$. We recall that
$$
\Vert F^{(n)}(\zeta)\Vert\le \frac{n!r(\Gamma)}{d(\Gamma,\Gamma_0)^{n+1}}\sup_{\theta\in\Gamma}\Vert F(\theta)\Vert,\,\,
\zeta\in\Delta_0,\,\,n\ge 0,
$$
where $r(\Gamma)$ is the radius of $\Gamma$, and $d(\Gamma,\Gamma_0)=\inf\{\vert\theta-\theta_0\vert;\theta\in\Gamma,\theta_0\in\Gamma_0\}$. Of course, these are {\it Cauchy's derivative inequalities}, with explicit constants, which
will be used in the sequel.  
\end{Rem}

% Lemma 4

\begin{Lem}\label{Cauchy_ineg} $(1)$ Let $D_+,\Delta_{0+},\Delta_+$ be disks in a conjugate symmetric open subset $U\subset\C$, such that 
$$D_+\subset\bar{D}_+\subset\Delta_{0+}\subset\bar{\Delta}_{0+}\subset\Delta_+\subset
\bar{\Delta}_+\subset\{\zeta\in\C;\Im{\zeta}>0\}.
$$
Let also $F\in\mathcal{O}_s(U,\M_2)$.  Then we have the estimates
$$
\Vert F^{(n)}_\H(q)\Vert\le  \left(\frac{2n!r(\Gamma_{0+})}{d(\Gamma_+,\Gamma_{0+})^{n+1}
d(\Gamma_{0+},C_{0+})}\right)\sup_{\zeta\in\Gamma_{+}}\Vert F(\zeta)\Vert,\,\,n\ge0,
$$
for all $q\in U_\H$ with $s_+(q)\in D_+$, where 
$C_{0+}, \Gamma_{0+},\Gamma_+$ are the boundaries of $D_+,
\Delta_{0+}, \Delta_+$, respectively.

$(2)$ Let $D,\Delta_{0},\Delta$ be disks in a conjugate symmetric open subset $U\subset\C,\break U\cap\R\neq\emptyset$, such that the centers of $D,\Delta_{0},\Delta$ belong to $\R$, and 
$$D\subset\bar{D}\subset\Delta_{0}\subset\bar{\Delta}_{0}\subset\Delta\subset\bar{\Delta}\subset U. 
$$
Let also $F\in\mathcal{O}_s(U,\M_2)$.  Then we have the estimates
$$
\Vert F^{(n)}_\H(q)\Vert\le \left(\frac{n!r(\Gamma_{0})}{d(\Gamma,\Gamma_{0})^{n+1}
d(\Gamma_{0},C)}\right)\sup_{\zeta\in\Gamma}\Vert F(\zeta)\Vert,\,\,n\ge0,
$$
for all $q\in U_\H$ with $\sigma(q)\subset D$, where 
$C, \Gamma_{0},\Gamma$ are the boundaries of $D,
\Delta_{0}, \Delta$, respectively.

\end{Lem}

{\it Proof.}\, (1) First of all,
for each subset $A_+\subset\{\zeta\in\C;\Im{\zeta}>0\}$, we 
put $A_-:=\{\bar{\zeta};\zeta\in A_+\}$.

Secondly, for every $q\in\H$ with $\sigma(q)\subset D_+\cup D_-$ 
the eigenvalues $s_\pm(q)$ of the normal operator $q$ on
$\C^2$ are distinct and we have a direct sum decomposition 
$\C^2=\C^2_+\oplus\C^2_-$, where $\C^2_{\pm}=\{\z\in\C^2;
q\z=s_{\pm}(q)\z\}$. Let also $E_\pm(q)$ be the orthogonal projection of $\C^2$ onto $\C^2_{\pm}$. Therefore,
$\Vert(\zeta{\bf I}-q)^{-1}E_\pm(q)\Vert=
\vert\zeta-s_\pm(q)\vert^{-1}$, whenever $\zeta\notin\sigma(q)$
(see also Remark \ref{calc_func}). 

Next, using Cauchy's inequalities, we obtain
$$
\left\Vert\frac{1}{2\pi i}\int_{\Gamma_{0+}} F^{(n)}(\zeta)(\zeta{\bf I}-q)^{-1}E_+(q)d\zeta\right\Vert\le
$$
$$
\left(\frac{n!r(\Gamma_{0+})}{d(\Gamma_+,\Gamma_{0+})^{n+1}}\right)\sup_{\zeta\in\Gamma_{+}}\Vert F(\zeta)\Vert\sup_{\zeta\in\Gamma_{0+}}\vert\zeta-s_+(q)\vert^{-1}\le 
$$
$$
\left(\frac{n!r(\Gamma_{0+})}{d(\Gamma_+,\Gamma_{0+})^{n+1}
d(\Gamma_{0+},C_{0+})}\right)\sup_{\zeta\in\Gamma_{+}}\Vert F(\zeta)\Vert,
$$
because $\sup_{\zeta\in\Gamma_{0+}}\vert\zeta-s_+(q)\vert^{-1}\le d(\Gamma_{0+},C_{0+})^{-1}$ when $s_+(s)\in D_{0+}$, where 
$C_{0+}, \Gamma_{0+},\Gamma_+$ are the boundaries of $D_+,
\Delta_{0+}, \Delta_+$, respectively.

Similarly,

$$
\left\Vert\frac{1}{2\pi i}\int_{\Gamma_{0-}} F^{(n)}(\zeta)(\zeta{\bf I}-q)^{-1}E_-(q)d\zeta\right\Vert\le
$$
$$
\left(\frac{n!r(\Gamma_{0-})}{d(\Gamma_-,\Gamma_{0-})^{n+1}
d(\Gamma_{0-},C_{0-})}\right)\sup_{\zeta\in\Gamma_{-}}\Vert F(\zeta)\Vert,
$$
where 
$C_{0-}, \Gamma_{0-},\Gamma_-$ are the boundaries of $D_,
\Delta_{0-}, \Delta_-$, respectively.

Because of symmetric conjugation, and the properties of $F$, we may set 
$$
M:=\sup_{\zeta\in\Gamma_{+}}\Vert F(\zeta)\Vert=\sup_{\zeta\in\Gamma_{-}}\Vert F(\zeta)\Vert,\,\,r:=r(\Gamma_{0+})=r(\Gamma_{0-}),
$$
$$
d:=d(\Gamma_+,\Gamma_{0+})=d(\Gamma_-,\Gamma_{0-}),\,\,
d_0=d(\Gamma_{0+},C_{0+})=d(\Gamma_{0-},C_{0-}),
$$
and as we have, 

$$
F^{(n)}_\H(q)=\frac{1}{2\pi i}\int_{\Gamma_{0+}} F^{(n)}(\zeta)(\zeta{\bf I}-q)^{-1}E_+(q)d\zeta+$$
$$
\frac{1}{2\pi i}\int_{\Gamma_{0-}} F^{(n)}(\zeta)(\zeta{\bf I}-q)^{-1}E_-(q)d\zeta,
$$
we infer that 
$$
\Vert F^{(n)}_\H(q)\Vert\le\frac{2rn!M}{d_0d^{n+1}},
$$ 
which are the desired estimates. 

(2) The proof is similar to that of part (1). In fact, following the lines of the previous proof, we have 
$$
\Vert F^{(n)}_\H(q)\Vert\le\left(\frac{n!r(\Gamma_{0})}{d(\Gamma,\Gamma_{0})^{n+1}
d(\Gamma_{0},C_{0})}\right)\sup_{\zeta\in\Gamma}\Vert F(\zeta)\Vert,
$$ 
because $\Vert(\zeta{\bf I}-q)\Vert^{-1}\le d(\Gamma_{0},C_{0})^{-1}$ when $\sigma(q)\subset D_0$. 

% Prop 1

\begin{Pro} Let $U\subset\C$ be a conjugate symmetric open set
containing $0$, and let $F\in\mathcal{O}_s(U,\M_2)$. 
Then there exists an open disk $V\subset U$ of center $0$ such that $Q((\lambda,0))\in V_\H$ for all $\lambda\in V$ and  
$$
F(\lambda)=\sum_{n\ge0}\frac{F^{(n)}_\H(q)}{n!} 
(\lambda{\bf I}-q)^n,\,\,q\in V_\H,
$$
where the series is absolutely convergent in $\M_2$. 
\end{Pro}

{\it Proof.}\, We choose the disks $D,\Delta_{0},\Delta$ as in Lemma \ref{Cauchy_ineg}(2). In addition, we assume that all of
them are centered at $0$. 

According to  Lemma \ref{Cauchy_ineg}(2), we deduce that
$$
\limsup_{n\to\infty}\left(\frac{\Vert F^{(n)}_\H(q)\Vert}{n!}
\right)^{1/n}\le d(\Gamma,\Gamma_0)^{-1}<\infty, 
$$
whenever $\sigma(q)\subset D$. Therefore, the formal series 
$$
G(\lambda)=\sum_{n\ge0}\frac{F^{(n)}_\H(q)}{n!} 
(\lambda{\bf I}-q)^n
$$
is absolutely convergent if $\Vert\lambda{\bf I}-q\Vert<\delta$, where $\delta=d(\Gamma,\Gamma_0)$.

Because this estimate does not depend on the radius of $D$, say $\rho$,  which shall determine $\rho$ such that $V=D\subset\Delta_0$. 

We assume that $0<2\rho<\delta$. 
If $\lambda\in V$, then $\bar{\lambda}\in V$, so $Q((\lambda,0))\in V_\H$. Next, if $q\in V_\H$, as both $s_\pm(q)\in V$, we 
must have $\vert\lambda-s_\pm(q)\vert<2\rho<\delta$,
implying $\Vert\lambda{\bf I}-q\Vert<\delta$. 

Note that
$$
G(\lambda)=\sum_{n\ge0}\frac{1}{2\pi in!}\left(\int_{\Gamma_0} 
F^{(n)}(\zeta)(\zeta{\bf I}-q)^{-1}d\zeta\right)(\lambda{\bf I}-q)^n=
$$
$$
\sum_{n\ge0}\frac{1}{2\pi i}\int_{\Gamma_0}\left(\frac{1}{2\pi i}\int_{\Gamma}
\frac{F(\theta)}{(\theta-\zeta)^{n+1}}d\theta\right)(\zeta{\bf I}-q)^{-1}(\lambda{\bf I}-q)^nd\zeta = 
$$
$$
\frac{1}{2\pi i}\int_{\Gamma}F(\theta)\left(\frac{1}{2\pi i}\int_{\Gamma_0}
\sum_{n\ge0}(\theta-\zeta)^{-n-1}(\lambda{\bf I}-q)^n(\zeta{\bf I}-q)^{-1}d\zeta\right)d\theta.
$$
Because $\Vert\lambda{\bf I}-q\Vert<\delta\le\vert\theta-\zeta\vert$ for all $\theta\in\Gamma$ and $\zeta\in\Gamma_0$, the following series is convergent, and  
$$
\sum_{n\ge0}(\theta-\zeta)^{-n-1}(\lambda{\bf I}-q)^n=
((\theta-\zeta-\lambda){\bf I}+q)^{-1}.
$$
Hence 
$$
\frac{1}{2\pi i}\int_{\Gamma_0}
\sum_{n\ge0}(\theta-\zeta)^{-n-1}(\lambda{\bf I}-q)^n
(\zeta{\bf I}-q)^{-1}d\zeta=(\theta-\lambda)^{-1}{\bf I},
$$
which implies the equality $G(\lambda)=F(\lambda)$. 
\vskip1cm

% Rem 13

\begin{Rem}\label{zeros}\rm  (1) Let $U\subset\C$ be a conjugate symmetric open set and let $F\in\mathcal{O}_s(U,\M_2)$ be arbitrary. We can easily describe the zeros of $F_\H$. 
Indeed, as we have $F_\H(q)=F(s_+(q))\nu_+(q)+F(s_-(q))
\nu_-(q)$, we have $F_\H(q)=0$ if and only if $F(s_\pm(q))=0$.
Therefore, setting $\mathcal{Z}(F):=\{\lambda\in U;F(\lambda)=0\}$,
and, similarly, $\mathcal{Z}(F_\H):=\{q\in U_\H;F_\H(q)=0\}$,
we must have
$$
\mathcal{Z}(F_\H)=\{q\in U_\H;\sigma(q)\subset\mathcal{Z}(F)\}. 
$$
In particular, if $U$ is connected and $\mathcal{Z}(F)$ has 
an accumulation point in $U$, then $F_\H=0$. 

(2) Let us observe that if $F\in\mathcal{O}_s(U,\M_2)$ has the property that $F_\H(x+y{\bf J})=0$ for all $x+iy\in U$, then
$F=0$, and so $F_\H=0$. Indeed, with the notation from
Remark \ref{extension}, because we have 
$$
F(\zeta)=\left(\begin{array}{cc} f_1(\zeta) & f_2(\zeta) 
\\ -\overline{f_2(\bar{\zeta})} & 
\overline{f_1(\bar{\zeta})}  \end{array}\right),\,
F(q_\zeta)=\left(\begin{array}{cc} f_1(\zeta) & f_2(\bar{\zeta}) 
\\ -\overline{f_2(\bar{\zeta})} & \overline{f_1(\zeta)}\end{array}\right),\,\,\zeta\in U,
$$ 
the assertion is clear.
\end{Rem}

% Rem 14    

\begin{Rem}\label{deriv}\rm Theorem \ref{H_afc} and its consequences  suggest 
a definition for $\H$-valued "analytic functions`` as elements of the set $\mathcal{R}(\Omega,\H)$, where $\Omega$ is a 
spectrally saturated open subset of $\H$. 
Because the expression "analytic function`` is quite improper 
in this context, the elements of $\mathcal{R}(\Omega,\H)$ will be called {\it Q-regular functions} on $\Omega$. In fact, 
the functions from $\mathcal{R}(\Omega,\H)$ may be also regarded as {\it Cauchy transforms} of the (stem) functions from $\mathcal{O}_s(U,\M_2)$, with $U=\mathfrak{S}(\Omega)$. 

We recall that there exists a large literature dedicated to a concept of "slice regularity``, which is a form of holomorphy in the context of 
quaternions (see for instance \cite{CoSaSt} and works quoted 
within). Till the end of this section we shall try to clarify
the connection between these concepts, showing that they
coincide on spectrally saturated open sets. 

For $\M_2$-valued functions defined on subsets of $\H$, the concept of slice regularity (see \cite{CoSaSt}) is defined as follows. 

Let $\S$ be the unit sphere of purely imaginary quaternions
(see Example \ref{spec_can_rep}).  
Let also $\Omega\in\H$ be an open set, and let $F:\Omega\mapsto\M_2$ be a differentiable function. In the spirit of \cite{CoSaSt}, we say that $F$ is {\it (right) slice regular} in $\Omega$ if for all $\mathfrak{s}\in\S$,
$$
 \bar{\partial}_\mathfrak{s}F(x{\bf I}+y\mathfrak{s}):
 =\frac{1}{2}
\left(\frac{\partial}{\partial x}+R_\mathfrak{s}\frac{\partial}{\partial y}\right)F(x+y\mathfrak{s})=0,
$$
on the set $\Omega\cap(\R{\bf I}+\R\mathfrak{s})$,
where $R_\mathfrak{s}$ is the right multiplication of the elements of $\M_2$ by $\mathfrak{s}$. 

Note that, unlike in \cite{GeSt}, we use the right slice regularity rather than the left one because of our regard to $\H$ as an algebra of operators on $\C^2$.

Of course, we are mainly interested by slice regularity of 
$\H$-valued functions, but the concept is valid for $\M_2$-valued
functions and plays an important role in our discussion.
\end{Rem}

\begin{Exa}\label{exa_sr}\rm 
(1) The convergent series of the form $\sum_{k\ge 0}a_kq^k$, 
on a set 
$\{q\in\H;\Vert q\Vert<r\}$, with $a_k\in\H$ for all $k\ge0$,
are $\H$-valued slice regular on their domain of definition. 
In fact, if actually $a_k\in\M_2$, such functions are $\M_2$-valued right slice regular on their domain of definition.

(2) The matrix Cauchy kernel on the open set $\Omega\subset\H$, defined by
$$
\Omega\ni q\mapsto (\zeta{\bf I}-q)^{-1}\in\M_2,
$$
is slice regular on $\Omega\subset\H$, whenever 
$\zeta\notin\mathfrak{S}(\Omega)$.
Indeed, for $q=x+y\mathfrak{s}\in\Omega\cap(\R{\bf I}+\R\mathfrak{s})$, we can write
$$
\frac{\partial}{\partial x}((\zeta-x){\bf I}-y\mathfrak{s})^{-1}=
((\zeta-x){\bf I}-y\mathfrak{s})^{-2},
$$
$$
R_\mathfrak{s}\frac{\partial}{\partial y}((\zeta-x){\bf I}-y\mathfrak{s})^{-1}=
 ((\zeta-x){\bf I}-y\mathfrak{s})^{-1}\mathfrak{s}((\zeta-x){\bf I}-y\mathfrak{s})^{-1}\mathfrak{s}=
$$
$$
-((\zeta-x){\bf I}-y\mathfrak{s})^{-2},
$$
because $\zeta{\bf I}$,  $\mathfrak{s}$ and 
$((\zeta-x){\bf I}-y\mathfrak{s})^{-1}$ commute in $\M_2$. Therefore,
$$
\bar{\partial}_\mathfrak{s}((\zeta{\bf I}-q)^{-1})=\bar{\partial}_\mathfrak{s}(((\zeta-x){\bf I}-y\mathfrak{s})^{-1})=0.
$$
\end{Exa}  

 % Lemma 5 

\begin{Lem}\label{right_reg} Let $\Omega\subset\H$ be a spectrally saturated open set. Then every function  $F\in \mathcal{R}(\Omega,\H)$ is right slice regular on $\Omega$.
\end{Lem}

{\it Proof.}\, We fix a function $F\in \mathcal{R}(\Omega,\H)$,
and let $U=\mathfrak{S}(\Omega)$. We use the representation 
$$
F_\H(q)=\frac{1}{2\pi i}\int_\Gamma F(\zeta)(\zeta{\bf I}-q)^{-1}d\zeta,\,\,q\in\Omega,\, \sigma(q)\subset\Delta,
$$
where $\Gamma$ is boundary of a Cauchy domain $\Delta\subset\bar{\Delta}\subset U$. Because we have
$$
\bar{\partial}_\mathfrak{s}((\zeta{\bf I}-q)^{-1})=\bar{\partial}_\mathfrak{s}(((\zeta-x){\bf I}-y\mathfrak{s})^{-1})=0
$$
for $q=x+y\mathfrak{s}\in\Omega\cap(\R{\bf I}+\R\mathfrak{s})$
as in  Exercise \ref{exa_sr}(2), we infer that 
$$
\bar{\partial}_\mathfrak{s}(F_\H(q))=\frac{1}{2\pi i}\int_\Gamma F(\zeta)\bar{\partial}_\mathfrak{s}((\zeta{\bf I}-q)^{-1})
d\zeta=0,
$$
which implies the assertion.

% Lemma 6

\begin{Lem}\label{splitting} Let $U\subset\C$ be a conjugate symmetric open set, let  $U_{\bf J}=\{x{\bf I}+y{\bf J}\in\H;x+iy\in U,x,y\in\R\}$, and let $f:U_{\bf J}\mapsto\H$ be 
such that $\bar{\partial}_{\bf J}f(x{\bf I}+y{\bf J})=0$. Then there are two functions $g,h:U_{\bf J}\mapsto\C_{\bf J}$ such that $\bar{\partial}_{\bf J}g=0,\,\bar{\partial}_{\bf J}h=0$ in $U_{\bf J}$,
and $f=g+{\bf L}h$, where $\C_{\bf J}=\{x{\bf I}+y{\bf J}\in\H;
x+iy\in U,x,y\in\R\}$.  
\end{Lem}

{\it Proof.}\, We proceed as in Lemma 4.1.7 in \cite{CoSaSt}.
We write $f=f_0{\bf I}+f_1{\bf J}+f_2{\bf K}+f_3{\bf L}$, 
where  $f_0,f_1,f_2,f_3$ are $\R$-valued functions. Then
$$
2\bar{\partial}_{\bf J}f(x{\bf I}+y{\bf J})=
\left(\frac{\partial}{\partial x}+R_{\bf J}\frac{\partial}{\partial y}\right)f(x{\bf I}+y{\bf J})=
$$ 

$$
\left(\frac{\partial f_0}{\partial x}-\frac{\partial f_1}{\partial y}\right){\bf I}+\left(\frac{\partial f_1}{\partial x}+\frac{\partial f_0}{\partial y}\right){\bf J}+
{\bf L}\left(\left(\frac{\partial f_3}{\partial x}-\frac{\partial f_2}{\partial y}\right){\bf I}+\left(\frac{\partial f_3}{\partial y}+\frac{\partial f_2}{\partial x}\right){\bf J}\right)=
$$ 
$$
2\bar{\partial}_{\bf J}(f_0{\bf I}+f_1{\bf J})+2{\bf L}(\bar{\partial}_{\bf J}(f_3{\bf I}+f_2{\bf J})=0. 
$$
Therefore, we may take $g=f_0{\bf I}+f_1{\bf J}$ and $h=
f_3{\bf I}+f_2{\bf J}$. 
\medskip

As mentioned in Remark \ref{consym}(5), a fixed conjugate symmetric open set  $U\subset\C$ can be associated with an axially symmetric set (see Definition 4.3.1 and Lemma 4.3.8 from
\cite{CoSaSt}), given by the formula
$$
\tilde{U}:=\{x{\bf I}+y\mathfrak{s};x+iy\in U,\,\mathfrak{s}\in\mathbb{S}\},
$$
which is the {\it circularization of} $U$ (as in \cite{GhMoPe}, Section 1.1).  

\begin{Pro}\label{UH_tildeU} For every  conjugate symmetric open set  $U\subset\C$ we have the equality  $U_\H=\tilde{U}.$
\end{Pro}

{\it Proof.}  If $q\in U_\H$, we can write $q=Q(\z^+(u))$
or $q=Q(\z^-(u))$, 
where $\z^\pm(u)=(x\pm i\sqrt{y^2-\vert u\vert^2}, u)\in\C^2$, and $x\pm iy\in U$,  for some complex number $u$ with $\vert u\vert\le\vert y\vert$, by Remark \ref{consym}.  As we have 
$$
Q(\z^\pm(u))=x{\bf I}\pm\sqrt{y^2-\vert u\vert^2}{\bf J}+u_1{\bf K}+u_2{\bf L}, 
$$
with $u=u_1+iu_2, \,u_1,u_2\in\R$, and for $y\neq 0$, we have
$$
\mathfrak{s}_\pm:=y^{-1}(\pm\sqrt{y^2-\vert u\vert^2}{\bf J}+u_1{\bf K}+u_2{\bf L})\in\mathbb{S},
$$
it follows that $q\in\tilde{U}$.

When $y=0$, then $x\in U\subset\tilde{U}$.
 
Conversely, 
let $q\in\tilde{U}$, so $q=x{\bf I}+y\mathfrak{s}$ for some 
$\mathfrak{s}\in\mathbb{S}$, and $x+iy\in U$. Of course,
$\mathfrak{s}=
a_1{\bf J}+
a_2{\bf K}+a_3{\bf L};a_1,a_2,a_3\in\R,a_1^2+a_2^2+a_3^2=1$.

To have
$q\in U_\H$, we must solve the equation,
$q=Q((x_1\pm i\sqrt{y_1^2-\vert u\vert^2},u)$, for some $x_1+iy_1\in U$, and $\vert u\vert^2\le y_1^2$, whose spectrum
is $\{x_1\pm iy_1\}$. On the other hand,
according to Example \ref{spec_can_rep}, the spectrum of $q$
is the set $\{x\pm iy\}$. 
Hence, we have the necessary conditions $x_1=x$ and $y_1=\pm y$. 
Note that we must have $u=a_2y+ia_3y$, and $a_1y=\pm\vert a_1y\vert$, which leads to a solution of the given equation for 
a suitable choice from $\{\pm y\}$. Consequently, $q\in U_\H$.

% Lemma 7 

\begin{Lem}\label{equiv-ons-dom} Let $U\subset\H$ be a conjugate symmetric open set, and let $\Phi_{\bf J}:U_{\bf J}\mapsto\H$ be  such that $\bar{\partial}_{\bf J}\Phi_{\bf J}=0$. Then there exists 
a function $\Phi\in\mathcal{R}(U_\H,\H)$ with
$\Phi_{\bf J}=\Phi\vert U_{\bf J}$.  
\end{Lem}

{\it Proof.}\, According to Lemma \ref{splitting}, we can write $\Phi_{\bf J}=
F_{\bf J}+{\bf L}G_{\bf J}$, with $F_{\bf J},\,G_{\bf J}:
U_{\bf J}\mapsto\C_{\bf J}$, and 
$\bar{\partial}_{\bf J}F_{\bf J}=0,\,\bar{\partial}_{\bf J}G_{\bf J}=0$ in $U_{\bf J}$. Note that we can write
$$
\C_{\bf J}=\left\lbrace\left(\begin{array}{cc} x+iy & 0 \\ 0 & x-iy
\end{array}\right); x,y\in\R\right\rbrace, 
$$ 
and  
$$
\bar{\partial}_{\bf J}=\frac{1}{2} \left(\begin{array}{cc}
\frac{\partial}{\partial x}+i\frac{\partial}{\partial y} & 0 \\
0 & \frac{\partial}{\partial x}-i\frac{\partial}{\partial y},
\end{array}\right).
$$

Because $\bar{\partial}_{\bf J}F_{\bf J}=0$, showing that the 
function $F_{\bf J}$ is analytic in $U_{\bf J}$, we have a local convergent series representation of this function under the form
$$
F_{\bf J}(x{\bf I}+y{\bf J})=\sum_{k\ge0}A_k((x-x_0)+(y-y_0){\bf J})^k
$$
in a neighborhood of each fixed point $x_0+y_0{\bf J}\in U_{\bf J}$, where $A_k\in\C_{\bf J}$ for all $k\ge0$. Using this local 
representation written in a matricial form, we derive the existence of an analytic function
$f_{\bf J}\in\mathcal{O}(U)$ such that
$$
F_{\bf J}(x{\bf I}+y{\bf J})=\left(\begin{array}{cc} f_{\bf J}(x+iy) & 0 \\ 0 & 
\overline{f_{\bf J}(x+iy)}
\end{array}\right),\,\, x+y{\bf J}\in U_{\bf J}.
$$

Similarly,
$$
G_{\bf J}(x{\bf I}+y{\bf J})=\left(\begin{array}{cc} g_{\bf J}(x+iy) & 0 \\ 0 & 
\overline{g_{\bf J}(x+iy)}\end{array}\right),\,\, x+y{\bf J}\in U_{\bf J},
$$
with $ g_{\bf J}\in\mathcal{O}(U)$. 
If 
$$
F(\zeta)=\left(\begin{array}{cc} f_{\bf J}(\zeta) & 0 \\ 0 & 
\overline{f_{\bf J}(\bar{\zeta})}\end{array}\right),\,
G(\zeta)=\left(\begin{array}{cc} g_{\bf J}(\zeta) & 0 \\ 0 & 
\overline{g_{\bf J}(\bar{\zeta})}\end{array}\right),\,\,
\zeta\in U,
$$
we have 
$F,G\in\mathcal{O}_s(U,\M_2)$, and $F_\H(q_\zeta)=
f_{\bf J}(x{\bf I}+y{\bf J}),\,G_\H(q_\zeta)=
g_{\bf J}(x{\bf I}+y{\bf J})$, with $\zeta=x+iy$,
via Remark \ref{extension}. Moreover, $F+{\bf L}G\in\mathcal{O}_s(U,\M_2)$, so setting  $\Phi(q)=F_\H(q)+{\bf L}G_\H(q)$ for all $q\in U_\H$, we have $\Phi\in\mathcal{R}(U_\H,\H)$ via Theorem \ref{H_afc}.

\begin{Thm}\label{equiv-ons-dom1} Let $\Omega\subset\H$ be a spectrally saturated open set, and let $\Phi:\Omega\mapsto\H$.
 The following conditions are equivalent:

$(i)$ $\Phi$ is  a slice regular function;  

$(ii)$ $\Phi\in\mathcal{R}(\Omega,\H)$, that is, $\Phi$ is 
$Q$-regular. 
\end{Thm}

{\it Proof.}\, If $\Phi\in\mathcal{R}(\Omega,\H)$, then $\Phi$
is slice regular, by Lemma \ref{right_reg}, so $(ii)\Rightarrow(i)$.

Conversely, let $\Phi$ be slice regular in $\Omega$. Then we 
have $\bar{\partial}_{\bf J}\Phi_{\bf J}=0$, where $\Phi_{\bf J}=\Phi\vert U_{\bf J}$. It follows from Lemma \ref{equiv-ons-dom}
that there exists $\Psi\in\mathcal{R}(U_\H,\H)$ with
$\Psi_{\bf J}=\Phi_{\bf J}$. This implies that 
$\Phi=\Psi$, because both $\Phi,\Psi$ are uniquely determined by 
$\Phi_{\bf J}, \Psi_{\bf J}$, respectively, the former by (the right side version of) Lemma 4.3.8 in\cite{CoSaSt}, and the latter by  Remark \ref{zeros}(2). Consequently, we also have $(i)\Rightarrow(ii)$.

\begin{Rem}\label{Martin}\rm 
An important integral formula extending Cauchy's formula 
to analytic functions in several variables is {\it Martinelli's formula}. Its first version appears in \cite{Mar},
and it was later independently obtained in \cite{Boc}. With our notation and framework, it can be stated in the following way. 

Let $D\subset\C^2$ be a bounded domain with piecewise-smooth
boundary $\partial D$, and lat $f:\bar{D}\mapsto\C$ be a function analytic in $D$ and continuous on $\bar{D}$. Then
$$
f({\bf w})=\int_{\partial D}f({\bf z})M(\bf{z,w})
$$
for all ${\bf w}\in\C^2$, where 
$$
M({\bf z,w})=\frac{1}{(2\pi i)^2}
\Vert Q({\bf z})-Q({\bf w})\Vert^{-2}(({\bar z}_1-{\bar w}_1)
d{\bar z}_2-({\bar z}_2-{\bar w}_2)
d{\bar z}_1)\wedge d{\bf z},
$$
with  $d{\bf z}:=dz_1\wedge dz_2$. 

In the last section we present a version of this formula, valid 
for commuting pairs of real operators.
\end{Rem}

\section{Quaternionic Spectrum of Real Operators}\label{QSRO}

For a real or complex Banach space $\mathcal{V}$, we denote by
$\mathcal{B(V)}$ the algebra of all bounded $\R$-(
respectively $\C$-)linear operators on $\mathcal{V}$.
If necessary, the identity on $\mathcal{V}$
will be denoted by ${\bf I}_{\mathcal{V}}$. Nevertheless, a multiple
of the identity $\lambda{\bf I}_{\mathcal{V}}$ will be often identified with the scalar $\lambda$, when no confusion is possible.  

Inspired by Definition 4.8.1 from \cite{CoSaSt}
(see also \cite{CoGaKi}; in fact, a similar idea goes back to Kaplansky, see \cite{Kul}), we give the folowing. 

% Def 3

\begin{Def}\label{Q-spectrum}\rm Let $\mathcal{X}$ be a real Banach space. For a given operator $T\in\mathcal{B(X)}$, the set 
$$
\rho_\H(T):=\{Q(\z)\in\H; \z=(z_1,z_2), (T^2-(z_1+\bar{z}_1)T+\vert z_1\vert^2+\vert z_2\vert^2)^{-1}\in\mathcal{B(X)}\}
$$
is said to be the {\it quaternionic resolvent} (or simply the
$Q$-{\it resolvent}) of $T$. 

The complement $\sigma_\H(T)=\H\setminus\rho_\H(T)$ is called 
the {\it quaternionic spectrum} (or simply the $Q$-{\it spectrum}) of $T$.
\end{Def}

We note that if $q=Q(\z)$, with $\z=(z_1,z_2)\in\C$, setting $\Re q=\Re z_1$, we have 
$$
T^2-(z_1+\bar{z}_1)T+\vert z_1\vert^2+\vert z_2\vert^2=T^2- 2\Re q\,T+\Vert q\Vert^2,
$$
and the right hand side is precisely the expression used in 
Definition 4.8.1 from \cite{CoSaSt}. Note also that Definition
\ref{Q-spectrum} applies to the class of $\R$-linear operators, which is larger that the class of $\H$-linear operators. 

% Rem 16

\begin{Rem}\label{ssspec}\rm
Looking at  Definition \ref{Q-spectrum}, we observe that if
$\{s_\pm(\bf z)\}$ are the eigenvaluse of $Q({\bf z})$, we have 
$Q(\z)\in\rho_\H(T)$ if and only if the operator
$$
T^2-(s_+({\bf z})+s_-({\bf z}))T+s_+({\bf z})s_-({\bf z})
$$
is invertible. Since this property depends only on the eigenvalues
of $Q({\bf z})$, it follows that $Q({\bf w})\in\rho_\H(T)$
whenever for  $Q(\z)\in\rho_\H(T)$ we have $\sigma(Q({\bf w}))=\sigma(Q({\bf z}))$, because this equality is equivalent to
$z_1+\bar{z}_1=w_1+\bar{w}_1$ and $\vert z_1\vert^2+\vert z_2\vert^2=
\vert w_1\vert^2+\vert w_2\vert^2$.
In particular, $Q(\z)\in\rho_\H(T)$ if and only if $Q(\z^*)\in\rho_\H(T)$, where $\z^*=(\bar{z}_1,-z_2)$ if $\z=(z_1,z_2)$ (see the second section), and $Q(\z)\in\sigma_\H(T)$ implies that
$Q((s_\pm(\z),0))\in\sigma_\H(T)$. In fact both sets $\rho_\H(T)$and $\sigma_\H(T)$ are spectrally saturated in $\H$ (see Remark \ref{consym}(2)).

The {\it complex spectrum} of operator $T\in\mathcal{B(X)}$ on the 
real Banach space $\mathcal{X}$ is given by 
$$
\sigma_\C(T):=\{\lambda\in\C;Q((\lambda,0))\in\sigma_\H(T)\}.
$$

Because the $Q$-spectrum of $T$ is spectrally saturated, we have
$\sigma_\H(T)=\sigma_\C(T)_\H$. We also have  $\lambda\in\sigma_\C(T)$ if and only if $\bar{\lambda}\in\sigma_\C(T)$. In addition 
$\lambda\in\rho_\C(T):=\C\setminus\sigma_\C(T)$ if and only if the operator $T^2-2\Re\lambda T +\vert\lambda\vert^2$ is invertible. 
\end{Rem}
\medskip

Let $\mathcal{X}$ be a real Banach space, and let $T\in\mathcal{B(X)}$. We denote by
$\mathcal{X}_\C$ the complexification of $\mathcal{X}$, written
as $\mathcal{X}_\C=\mathcal{X}\oplus i\mathcal{X}$, or simply as $\mathcal{X}+i\mathcal{X}$. The operator
$T$ can be extended to $\mathcal{X}_\C$ via the formula
$T_\C(x+iy)=Tx+iTy$ for all $x,y\in\mathcal{X}$. It is clear that $T_ \C$ is a bounded $\C$-linear operator.  Let
 $T_\C^{(2)}$ be the $2\times2$ diagonal operator with $T_\C$ on the diagonal, acting on  $\mathcal{X}_\C^2:=
\mathcal{X}_\C\oplus\mathcal{X}_\C$. As for every $\z\in\C^2$ the matrix $Q(\z)$ also acts on  $\mathcal{X}_\C^2$, we may state the
following.

% Lemma 8 

\begin{Lem}\label{Q-sp} A quaternion $Q(\z)$ is in the set $\rho_\H(T)$ if and only if the operators $T_\C^{(2)}-Q(\z)$
and $T_\C^{(2)}-Q(\z^*)$ are invertible in 
$\mathcal{B({X}}_\C^2)$.
\end{Lem}

{\it Proof}. The assertion follows via the equalities
$$
\left(\begin{array}{cc} T_\C-z_1 & -z_2 \\
\bar{z}_2 & T_\C-\bar{z}_1\end{array}\right)
\left(\begin{array}{cc} T_\C-\bar{z}_1 & z_2 \\
-\bar{z}_2 & T_\C-z_1\end{array}\right)=
$$
$$
\left(\begin{array}{cc} T_\C-\bar{z}_1 & z_2 \\
-\bar{z}_2 & T_\C-z_1\end{array}\right)
\left(\begin{array}{cc} T_\C-z_1 & -z_2 \\
\bar{z}_2 & T_\C-\bar{z}_1\end{array}\right)=
$$
$$
\left(\begin{array}{cc} (T_\C-z_1)(T_\C-\bar{z}_1)+\vert z_2\vert^2 & 0 \\
0  & \vert z_2\vert^2+(T_\C-\bar{z}_1)(T_\C-z_1)\end{array}\right).
$$

Consequently, the operators $T_\C^{(2)}-Q(\z),\,T_\C^{(2)}-
Q(\z^*)$ are invertible in
$\mathcal{B}({\mathcal X}_\C^2)$ if and only if the operator 
$$(T_\C-z_1)(T_\C-\bar{z}_1)+\vert z_2\vert^2
=T_\C^2-(z_1+\bar{z}_1)T_\C+\vert z_1\vert^2+\vert z_2\vert^2
$$ 
is invertible in  $\mathcal{B}(\mathcal{X}_\C)$.  Because we have
$$
T_\C^2-(z_1+\bar{z}_1)T_\C+\vert z_1\vert^2+\vert z_2\vert^2=
(T^2-(z_1+\bar{z}_1)T+\vert z_1\vert^2+\vert z_2\vert^2)_\C,
$$
the operators $T_\C^{(2)}-Q(\z),\,T_\C^{(2)}-
Q(\z^*)$ are invertible in
$\mathcal{B}({\mathcal X}_\C^2)$ if and only if the operator $T^2-(z_1+\bar{z}_1)T+\vert z_1\vert^2+\vert z_2\vert^2$
is invertible in $\mathcal{B(X)}$.

% Ex 3

\begin{Exa}\rm One of the simplest possible example is to take 
$\mathcal{X}=\R$ and $T$ the operator $Tx=\tau x$ for
all $x\in\R$, where $\tau\in\R$ is fixed. We have $\mathcal{X}_\mathbb{C}=\C$, and $T_\C$ is given by the same formula, acting 
on $\C$.

The $Q$-spectrum
of $T$ is then the set 
$$\{Q(\z);\z=(z_1,z_2)\in\C, 
\tau^2-\tau(z_1+\bar{z}_1)+\vert z_1\vert^2+\vert z_2\vert^2
=0\}=
$$
$$
\{Q(\z);\z=(z_1,z_2)\in\C,\Re z_1=\tau, \Im z_1=z_2=0 \}=
\{Q((\tau,0))\}. 
$$
Consequently, $\sigma_\H(T)=\{Q((\tau,0))\}$, and $\sigma_\C(T)=\{\tau\}$. 

More examples can be found in \cite{CoSaSt}, Section 4.9, which can be
easily adapted to our context. For instance (see Example 4.9.3 from
\cite{CoSaSt}), if an operator $M$  is given by the matrix
$$
M=\left(\begin{array}{cc} a & 0 \\ 0 & b \end{array}\right),\,\,a,b\in\R,
$$
acting on $\mathcal{X}=\R^2$, we identify $\mathcal{X}_\C$ with $\C^2$,
and $M_\C$ is given by the same matrix, acting on $\C^2$. We have $Q(\z)\in\sigma_\H(M),\,\z=(z_1,z_2),$
if and only if $\z$ is a solution of at least one of the equations
$$
a^2-2\Re z_1\,a+\vert z_1\vert^2+\vert z_2\vert^2=0;\,\,
b^2-2\Re z_1\,b+\vert z_1\vert^2+\vert z_2\vert^2=0.
$$
Therefore, we should have either $\Re z_1=a,\,\Im z_1=z_2=0$, or 
$\Re z_1=b,\,\Im z_1=z_2=0$. Hence $\sigma_\H(M)=\{Q(a,0)),Q(b,0))\}.$
Moreover, $\sigma_\C(M)=\{a,b\}$. 

Note that $\sigma_\C(T)=\sigma(T_\C)$, and $\sigma_\C(M)=\sigma(M_\C)$.
This is a general property, proved in Lemma \ref{specr_specc}.
\end{Exa}

% Rem 17

\begin{Rem}\label{conjug}\rm Let again $\mathcal{X}$ be a real Banach space, and let 
$\mathcal{X}_\C$ be its complexification. For every $u=x+iy\in\mathcal{X}_\C$ with $x,y\in\mathcal{X}$ we put $\bar{u}=
x-iy$. In other words, the map $\mathcal{X}_\C\ni u\mapsto\bar{u}\in\mathcal{X}_\C$ is a conjugation, also denoted
by $C$. Hence $C$ is $\R$-linear and $C^2$ is the identity on
$\mathcal{X}_\C$. 

Fixing an operator $S\in\mathcal{B}(\mathcal{X}_\C)$, we define 
the operator $S^\flat\in\mathcal{B}(\mathcal{X}_\C)$ to be equal to $CSC$. It is easily seen that the map 
$\mathcal{B}(\mathcal{X}_\C)\ni S\mapsto S^\flat\in \mathcal{B}(\mathcal{X}_\C)$ is a unital conjugate-linear automorphism, whose square is the identity on $\mathcal{B}(\mathcal{X}_\C)$. 
Because $\mathcal{X}=\{u\in\mathcal{X}_\C; Cu=u\}$, we have
$S^\flat=S$ if and only if $S(\mathcal{X})\subset\mathcal{X}$.
In particular, $\mathcal{X}$ is invariant under $S+S^\flat$,
and if $T\in\mathcal{B}(\mathcal{X})$, we have $T_\C^\flat=
T_\C$. 
\end{Rem}

% Lem 9

\begin{Lem}\label{specr_specc} Let $\mathcal{X}$ be a real Banach space, and 
let $T\in\mathcal{B(X)}$. We have the equality $\sigma_\C(T)=
\sigma(T_\C)$. 
\end{Lem}

{\it Proof.}\, We first recall that $\lambda\in\sigma(T_\C)$
if and only if $\bar{\lambda}\in\sigma(T_\C)$ (see Remark \ref{ssspec}). Next, if $\lambda\notin\sigma(T_\C)$, then
$$
(T_\C^{(2)}-Q((\lambda,0)))^{-1}=\left(\begin{array}{cc}
(T_\C-\lambda)^{-1} & 0 \\ 0 & (T_\C-\bar{\lambda})^{-1}
\end{array}\right),
$$
showing that $\lambda\notin\sigma_\C(T)$. Hence,
$\sigma_\C(T)\subset\sigma(T_\C)$.

Conversely, if $\lambda\notin\sigma_\C(T)$, so 
$\bar{\lambda}\notin\sigma_\C(T)$, we have  
$Q((\lambda,0)),Q((\bar{\lambda},0)\notin\sigma_\H(T)$. 
Then both $(T_\C^{(2)}-Q((\lambda,0)))^{-1},(T_\C^{(2)}-
Q((\bar{\lambda},0)))^{-1}$ belong to
$\mathcal{B({X}}_\C^2)$. Hence, as in Lemma \ref{Q-sp}, the operator
$$
T_\C^2-2\Re\lambda T_\C+\vert\lambda\vert^2=(T_\C-\lambda)
(T_\C-\bar{\lambda})
$$
must be invertible, implying that $T_\C-\lambda$ is invertible, so $\lambda\notin \sigma(T_\C)$.

% Rem 18

\begin{Rem}\label{comp_spec}\rm If $\mathcal{X}$ is a real Banach space, and $T\in\mathcal{B(X)}$, both $\sigma_\H(T)$
and $\sigma_\C(T)$ are nonempty compact subset of 
$\H,\,\C$, respectively. Indeed, if $\Vert Q(\z)\Vert=\Vert Q(\z^*)\Vert>\Vert T_\C^{(2)}\Vert$, both operators 
$T_\C^{(2)}-Q(\z)$ and $T_\C^{(2)}-Q(\z^*)$ are invertible in 
$\mathcal{B({X}}_\C^2)$, via a classical argument of perturbation theory (see \cite{DuSc}, VII.6.1). Therefore, the set $\sigma_\H(T)$ is bounded, because $T_\C^{(2)}$ is bounded. 
Next, if $Q(\z),\,Q(\z^*)\in\rho_\H(T)$, the operator
$T_\C^2-(z_1+\bar{z}_1)T_\C+\vert z_1\vert^2+\vert z_2\vert^2$
must be invertible, as in the proof of Lemma \ref{Q-sp}. 
Therefore, there is an $\epsilon>0$ such that $\Vert Q({\bf w})
-Q(\z)\Vert<\epsilon$, with ${\bf w}=(w_1,w_2)$, implying that the 
operator $T_\C^2-(w_1+\bar{w}_1)T_\C+\vert w_1\vert^2+\vert w_2\vert^2$ is invertible too, by continuity reasons.  Hence $\rho_\H(T)$ is an
open set (via Lemma \ref{Q-sp}), showing that $\sigma_\H(T)$
is compact. Finally, since $\sigma_\C(T)=\sigma(T_\C)$, and the 
latter is compact and nonempty, it follows that $\sigma_\H(T)$
is also nonempty. 
\end{Rem}

% Rem 19

\begin{Rem}\label{afcro}\rm If $\mathcal{X}$ is a real Banach space, and 
$T\in\mathcal{B(X)}$, we have the usual analytic functional calculus for the operator $T_\C\in\mathcal{B}(\mathcal{X}_\C)$.
That is, if $U\supset\sigma(T_\C)$ is an open set in $\C$ and 
$F:U\mapsto\mathcal{B}(\mathcal{X}_\C)$ is analytic, we may put
$$
F(T_\C)=\frac{1}{2\pi i}\int_\Gamma F(\zeta)(\zeta-T_\C)^{-1}
d\zeta,
$$  
where $\Gamma$ is the boundary of a Cauchy domain containing  
$\sigma(T_\C)$ in $U$. In fact, because $\sigma(T_\C)$ is conjugate symmetric, we may and shall assume 
that both $U$ and $\Gamma$ are conjugate symmetric. A natural 
question is when we have $F(T_\C)^\flat=F(T_\C)$, which would
imply the invariance of $\mathcal{X}$ under $F(T_\C)$. 
\end{Rem}

With the conditions of the remark from above we have the following.

% Thm 4

\begin{Thm}\label{afcro1} If $F:U\mapsto\mathcal{B}(\mathcal{X}_\C)$ is analytic and $F(\zeta)^\flat=F(\bar{\zeta})$ for all $\zeta\in U$, then  $F(T_\C)^\flat=F(T_\C)$ for all $T\in
\mathcal{B}(\mathcal{X})$.
\end{Thm}

{\it Proof.}\,  As in Remark \ref{afcro}, we have 
$$
F(T_\C)=\frac{1}{2\pi i}\int_\Gamma F(\zeta)(\zeta-T_\C)^{-1}
d\zeta,
$$  
where both $U$ and $\Gamma$ are conjugate symmetric. We put 
$\Gamma_\pm:=\Gamma\cap\C_\pm$, where $\C_+$ (resp. $\C_-$) equals to $\{\lambda\in\C;\Im\lambda\ge0\}$ (resp. $\{\lambda\in\C;\Im\lambda\le0\}$). We write $\Gamma_+=\cup_{j=1}^m\Gamma_{j+}$, where $\Gamma_{j+}$ are the connected components
of $\Gamma_+$. Similarly, we write  $\Gamma_-=\cup_{j=1}^m\Gamma_{j-}$, where $\Gamma_{j-}$ are the connected components
of $\Gamma_-$, and $\Gamma_{j-}$ is the reflexion of 
$\Gamma_{j+}$ with respect of the real axis. 

As $\Gamma$ is a finite union of Jordan piecewise smooth closed curves, 
for each index
$j$ we have a parametrization $\phi_j:[0,1]\mapsto\C$ such that
$\phi_j([0,1])=\Gamma_{+j}$. Taking into account the positive orientation, the function $t\mapsto-\overline{\phi_j(t)}$ is a 
parametrization of $\Gamma_{-j}$. Setting $\Gamma_j=\Gamma_{+j}\cup\Gamma_-$, we can write
$$
F_j(T_\C):=\frac{1}{2\pi i}\int_{\Gamma_j} F(\zeta)(\zeta-T_\C)^{-1}
d\zeta=
$$
$$
\frac{1}{2\pi i}\int_0^1 F(\phi_j(t))(\phi_j(t)-T_\C)^{-1}
\phi_j'(t)dt
$$
$$
-\frac{1}{2\pi i}\int_0^1 F(\overline{\phi_j(t)})(\overline{\phi_j(t)}-T_\C)^{-1}\overline{\phi_j'(t)}dt.
$$
Therefore,
$$
F_j(T_\C)^\flat=
-\frac{1}{2\pi i}\int_0^1 F(\phi_j(t))^\flat(\overline{\phi_j(t)}-T_\C)^{-1}\overline{\phi_j'(t)}dt
$$
$$
+\frac{1}{2\pi i}\int_0^1 F(\overline{\phi_j(t)})^\flat(\phi_j(t)-T_\C)^{-1}\phi_j'(t)dt.
$$
According to our assumption on the function $F$, we obtain
$F_j(T_\C)=F_j(T_\C)^\flat$ for all $j$, and therefore 
$$
F(T_\C)^\flat=\sum_{j=1}^mF_j(T_\C)^\flat=\sum_{j=1}^mF_j(T_\C)=F(T_\C) $$.

% Ex 4

\begin{Exa}\label{ex2}\rm Let $\mathcal{X}=\R^2$, so $\mathcal{X}_\C=\C^2$.
Let us first observe that  
$$
S=\left(\begin{array}{cc} a_1 & a_2 \\ a_3 &  a_4
\end{array}\right)\,\,\Longleftrightarrow
S^\flat=\left(\begin{array}{cc} \bar{a}_1 & \bar{a}_2 \\ 
\bar{a}_3 &  \bar{a}_4\end{array}\right),
$$
for all $a_1,a_2,a_3,a_4\in\C$. 

Next we consider the operator $T\in\mathcal{B}(\R^2)$ given by the matrix 
$$
T=\left(\begin{array}{cc} u & v \\ -v & u
\end{array}\right),
$$
where $u,v\in\R, v\neq0$. The extension $T_\C$ of the operator $T$ to $\C^2$ is given by the same formula. In fact, $T_\C$ is a quaternion,
denoted by $\tau$, for simplicity. Note that
$$
\sigma(\tau)=\{\lambda\in\C;(\lambda-u)^2+v^2=0\}=
\{u\pm iv\}.
$$
Let $U\subset\C$ be an open set with $U\supset\{u\pm iv\}$,
and let $F:U\mapsto\M_2=\mathcal{B}(\C^2)$ be analytic. 
We shall compute explicitly $F(\tau)$. Assuming $v>0$, we have 
$s_\pm(\tau)=u\pm iv$, and $\nu_\pm(\tau)=(\sqrt{2})^{-1}(1,\pm i)$. 
As in the proof of Lemma \ref{vect_fc}, we have
$$
F_\H(\tau)=F(s_+(\tau))E_+(\tau)+F(s_-(\tau))E_-(\tau)=
$$
$$
\frac{1}{2}F(u+iv)\left(\begin{array}{cc} 1 & -i \\ i & 1 \end{array}\right)+\frac{1}{2}F(u-iv)\left(\begin{array}{cc} 1 & i \\ -i & 1 \end{array}\right),
$$
because
$$
E_\pm(\tau){\bf w}=\langle {\bf w},\nu_\pm(\tau)\rangle \nu_\pm(\tau){\bf w}=\frac{1}{2}\left(\begin{array}{cc} 1 & \mp i \\ \pm i & 1 \end{array}\right)\left(\begin{array}{c} w_1 \\ w_2 \end{array}\right),
$$ 
for all ${\bf w}=(w_1,w_2)\in\C^2$.

Therefore,
$$
F(T_\C)=\frac{F(u+iv)}{2}\left(\begin{array}{cc} 1 & -i \\ i & 1
\end{array}\right)+\frac{F(u-iv)}{2}\left(\begin{array}{cc} 1 & i \\ -i & 1
\end{array}\right).
$$

Assume now that $F(T_\C)^\flat=F(T_\C)$. Then we must have

$$
(F(u+iv)-F(u-iv)^\flat)\left(\begin{array}{cc} 1 & -i \\ i & 1
\end{array}\right)= 
(F(u+iv)^\flat-F(u-iv))\left(\begin{array}{cc} 1 & i \\ -i & 1
\end{array}\right).
$$
We also need the equalities
$$
\left(\begin{array}{cc} 1 & -i \\ i & 1
\end{array}\right)\left(\begin{array}{c} 1 \\ i \end{array} \right)=2\left(\begin{array}{c} 1 \\ i \end{array}\right),\,\,
\left(\begin{array}{cc} 1 & -i \\ i & 1
\end{array}\right)\left(\begin{array}{c} 1 \\ -i \end{array} \right)=0,
$$
$$
\left(\begin{array}{cc} 1 & i \\ -i & 1
\end{array}\right)\left(\begin{array}{c} 1 \\ -i \end{array} \right)=2\left(\begin{array}{c} 1 \\ -i \end{array}\right),\,\,
\left(\begin{array}{cc} 1 & i \\ -i & 1
\end{array}\right)\left(\begin{array}{c} 1 \\ i \end{array} \right)=0,
$$
With this remark we finally deduce that
$$
(F(u+iv)-F(u-iv)^\flat)\left(\begin{array}{c} 1 \\ i \end{array}\right)=0,
$$
and
$$
(F(u-iv)-F(u+iv)^\flat)\left(\begin{array}{c} 1 \\ -i \end{array}\right)=0,
$$
which are necessary conditions for the equality $F(T_\C)^\flat=F(T_\C)$. 
As a matter of fact, this example shows that the condition used in Theorem \ref{afcro1} is sufficient but it might not be necessary. 
\end{Exa}

% Rem 20

\begin{Rem}\label{new_hol_sp}\rm  (1) Let $U\subset\C$ be a conjugate symmetric open set, and
let $\mathcal{X}$ be a real Banach space. We denote by 
$\mathcal{O}_c(U,\mathcal{B}(\mathcal{X}_\C))$ the set of all analytic maps
$F:U\mapsto\mathcal{B}(\mathcal{X}_\C)$ such that $F(\zeta)^\flat=F(\bar{\zeta})$ for all $\zeta\in U$. When $\mathcal{X}=\R$,
we put $\mathcal{O}_c(U,\mathcal{B}(\mathcal{X}_\C))=\mathcal{O}_c(U)$. In this case, we have $\mathcal{O}_c(U)=\mathcal{O}_s(U)$,
and $\mathcal{O}_c(U,\mathcal{B}(\mathcal{X}_\C))$ is a $\mathcal{O}_c(U)$-module.

Moreover, $\mathcal{O}_c(U,\mathcal{B}(\mathcal{X}_\C))$ is a unital $\R$-algebra, 
containing all polynomials $P(\zeta)=\sum_{k=0}^m (A_k)_\C\zeta^k$, with $A_k\in \mathcal{B}(\mathcal{X})$.

(2) The injective linear  map $\M_2\ni a\mapsto M_a\in\mathcal{B(X)}$, given $M_ab=ab,\,b\in\M_2$ induces an injective linear map of $\mathcal{O}_s(U,\M_2))$ into $\mathcal{O}_c(U,\mathcal{B(\M}_2))$. Specifically, given $F\in\mathcal{O}_s(U,\M_2))$, that is, an analytic stem functions (see Remark \ref{ssc_vs_sf}), we have $M_F\in \mathcal{O}_c(U,\mathcal{B(\M}_2))$, where $M_F(\zeta)b=F(\zeta)b$ for all
$\zeta\in U$ and $b\in\M_2$. Indeed, as we have $\M_2=\H+i\H$,
with the conjugation $C:u+iv\mapsto u-iv,\,u,v\in\H$, 
and writing $F(\zeta)=F_1(\zeta)+ iF_2(\zeta)$ with $F_1,F_2$
$\H$-valued, for $b=u+iv, u,v\in\H,$ we have
$$
M_{F(\zeta)}^\flat b=CM_{F(\zeta)}Cb=
C(F_1(\zeta)+iF_2(\zeta))(u-iv)=
$$
$$
C(F_1(\zeta)u+F_2(\zeta)v+i(F_2(\zeta)u-F_1(\zeta)v)=
$$
$$
(F_1(\zeta)-iF_2(\zeta))(u+iv)=M_{F(\bar{\zeta})}b,
$$
for every $\zeta\in U$. This remark shows that the space 
$\mathcal{O}_s(U,\M_2))$ may be regarded as a subspace of $\mathcal{O}_c(U,\mathcal{B(\M}_2))$. 

Fixing $F\in \mathcal{O}_c(U,\mathcal{B(X}_\C))$, we must have
%%%%%%%%%%%%%%%%%%%%%%%%%%%%%%%%%%%%%%%%%%%%%%%%%%%%%%%%%%%%%%%%
% \mathcal{X}_\C was replaced by \mathcal{B(X}_\C)
%%%%%%%%%%%%%%%%%%%%%%%%%%%%%%%%%%%%%%%%%%%%%%%%%%%%%%%%%%%%
$F(T_\C)^\flat=F(T_\C)$ for all $T\in\mathcal{B}(\mathcal{X})$. 
This allows us to define 
$F(T)=F(T_\C)\vert\mathcal{X}$, because $\mathcal{X}$ is invariant 
under $F(T_\C)$. In addition, we have the following.
\end{Rem}

% Thm 5

\begin{Thm}\label{afcro2} For every $T\in\mathcal{B}(\mathcal{X})$, the map
$$
\mathcal{O}_c(U,\mathcal{B(X}_\C))\ni F\mapsto F(T)\in\mathcal{B}(\mathcal{X})
$$
%%%%%%%%%%%%%%%%%%%%%%%%%%%%%%%%%%%%%%%%%%%%%%%%%%%%%%%%%%%%%%%%
% \mathcal{X}_\C was replaced by \mathcal{B(X}_\C)
%%%%%%%%%%%%%%%%%%%%%%%%%%%%%%%%%%%%%%%%%%%%%%%%%%%%%%%%%%%%
is $\R$-linear and has the property $(Ff)(T)=F(T)f(T)$ for all
$f\in\mathcal{O}_c(U)$ and $F\in \mathcal{O}_c(U,\mathcal{B}(\mathcal{X}_\C))$. 
Moreover,
$P(T)=\sum_{k=0}^m A_kT^k$ for any polynomial 
$P(\zeta)=\sum_{k=0}^m A_k\zeta^k,\,\zeta\in\C$, with coefficients $A_k\in \mathcal{B}(\mathcal{X})$.
\end{Thm}

{\it Proof.}\,The assertion is a direct consequence of the 
analytic functional calculus of the operator $T_\C$, and 
Theorem \ref{afcro1}. We omit the details.
\medskip

Remark \ref{new_hol_sp} shows  that functional calculus
given by Theorem \ref{H_afc} is compatible with that given by 
Theorem \ref{afcro2}.

% Def 3

\begin{Def} Let $\mathcal{X}$ be a real Banach space. We say that
$\mathcal{X}$ is a {\it (left) $\H$-module} if there exists 
a unital $\R$-algebra morphism of $\H$ into 
$\mathcal{B}(\mathcal{X})$. In this case, the elements of 
$\H$ in $\mathcal{B}(\mathcal{X})$ may regarded as $\R$-linear operators. 
\end{Def}

% Cor 6

\begin{Cor}\label{afcqo} If $\mathcal{X}$ is a $\H$-module, for every polynomial $P(\zeta)=\sum_{k=0}^m A_k\zeta^k,\,\zeta\in\C$, with coefficients $A_k\in\H$, we have 
$P(T)=\sum_{k=0}^m A_kT^k$.
\end{Cor}

% Ex 5 

\begin{Exa}\label{spec_lm}\rm
Given a measurable space $(\Omega,\Sigma)$, and a positive
measure $\mu$ on $\Sigma$, one of the most interesting quaternionic  spaces seems to be the space $L^2(\Omega,\mu;\H)$, 
consisting of $\Sigma$-measurable $\H$-valued functions $F$, with the property $\Vert F\Vert_2^2:=\break\int_\Omega\Vert F(\omega)\Vert^2d\mu(\omega)
<\infty,\,F\in L^2(\Omega,\mu;\H)$. In fact this space admits a $\H$-valued inner product given by
$$
\langle F,G\rangle_2=\int_\Omega F(\omega)G(\omega)^*d\mu(\omega),
\,\, F,G\in L^2(\Omega,\mu;\H),
$$
which induces the norm $\Vert *\Vert_2^2$. Of course, with respect to this 
norm, $L^2(\Omega,\mu;\H)$ is, in particular, an $\R$-Banach space.
It is also a bilateral $\H$-module.

Because $\M_2=\H+i\H$, if $\mathcal{X}=L^2(\Omega,\mu;\H)$, we may identify the space 
$\mathcal{X}_\C$ with the space $L^2(\Omega,\mu;\M_2)$, consisting of 
$\Sigma$-measurable $\M_2$-valued functions $F$, with the property $\int_\Omega\Vert F(\omega)\Vert^2d\mu(\omega)$,
which is the natural extension of the previous norm.

Next, let $\Theta:\Omega\mapsto\H$ be a $\mu$-essentially bounded
function. We may define on $L^2(\Omega,\mu;\H)$ the operator
$TF(\omega)=\Theta(\omega)F(\omega),\,\omega\in\Omega,\,F\in L^2(\Omega,\mu;\H)$, which is $\R$-linear. Note also that
$$
\Vert TF\Vert_2^2=\int_\Omega \Theta(\omega)F(\omega)F(\omega)^*\Theta(\omega)^*d\mu(\omega)=
$$
$$
\int_\Omega \Vert\Theta(\omega)\Vert^2\Vert F(\omega)\Vert^2 
d\mu(\omega)\le\Vert\Theta\Vert_\infty^2\int_\Omega\Vert F(\omega)\Vert^2 d\mu(\omega),
$$
where $\Vert\Theta\Vert_\infty$ is the essentiel upper bound 
of $\Vert\Theta(\omega)\Vert$ on $\Omega$, showing that the
operator $T$ is bounded.

The space $L^2(\Omega,\mu;\H)$ has also a right $\H$-module structure induced by the 
map $\H\ni q\mapsto R_q\in\mathcal{B}(L^2(\Omega,\mu;\H))$, given by 
$(R_qF)(\omega)=F(\omega)q,\,\omega\in\Omega$. With respect to this
structure, the operator $T$ is right $\H$-linear. 

The operator $T_\C$, that is the natural extension of $T$ to 
$L^2(\Omega,\mu;\M_2)$, is given by the same formula, written 
for functions $F\in L^2(\Omega,\mu;\M_2)$. 

Let us compute the $Q$-spectrum of $T$. According to Definition \ref{Q-spectrum}, we have
$$
\rho_\H(T)=\{q\in\H; (T^2-2\Re q\,T+\Vert q\Vert^2)^{-1}\in\mathcal{B}(\mathcal{X})\}.
$$
Consequently, $q\in\sigma_\H(T)$ if and only if zero is in the (essential) range of the function 
$$
\tau(q,\omega):=\Theta(\omega)^2-2\Re q\,\Theta(\omega)+\Vert q\Vert^2
,\,\omega\in\Omega.
$$ 
In other words, we must have
$$
\sigma_\H(T)=\{q\in\H;\forall\epsilon>0\,\,\mu(\{\omega;\Vert
\tau(q,\omega)\Vert<\epsilon\})>0\},
$$
and so
$$
\sigma_\C(T)=\{\zeta\in\C;\forall\epsilon>0\,\,\mu(\{\omega;\Vert
\tau(q_\zeta,\omega)\Vert<\epsilon\})>0\},
$$
where $q_\zeta=Q((\lambda,0)). $

In particular, if $P(\zeta)=\sum_{j=1}^m R_{q_j}\zeta^j$, we have 
$P(T)=\sum_{j=1}^m R_{q_j}T^j\in \mathcal{B}(L^2(\Omega,\mu;\H))$,
where $q_1,\ldots,q_m\in\H$. Of course,  for every open conjugate symmetric subset
$U\subset\C$ containing $\sigma_\C(T)$, and for every function 
$F\in\mathcal{O}_c(U,\mathcal{B}(\mathcal{X}_\C))$, we may 
construct the operator $F(T)\in\mathcal{B}(L^2(\Omega,\mu;\H))$.

Assuming  $\Omega$  a Hausdorff compact space and $\Theta$
continuous, the $Q$-spectrum of $T$ is given by the  set
$$
\sigma_\H(T)=\{q\in\H; \{q,q^*\}\cap\Theta(\Omega)\neq\emptyset\},
$$
via Lemma \ref{Q-sp}. In fact, it is the spectrum $\sigma_\C(T)$
which can be used for the computation of the analytic functional calculus of $T$. If we write $\Theta(\omega)$ as $Q((\theta_1(\omega),
\theta_2(\omega)),\,\omega\in\Omega,$ for some continuous functions 
$\theta_1,\theta_2:\Omega\mapsto\C$, we have the equality
$$
\sigma_\C(T)=\{\zeta\in\C;\exists \omega\in\Omega:
\zeta^2-2\Re(\theta_1(\omega))\zeta+\vert\theta_1(\omega)\vert^2+
\vert\theta_2(\omega)\vert^2=0\}.
$$
\end{Exa}

\section{Quaternionic Joint Spectrum of Paires}
\label{QJSP}

A strong connection between the spectral theory of pairs of 
commuting operators in a complex Hilbert space  and the algebra of quaternions has been firstly noticed in \cite{Vas1}. Other
connections will be presented in this section.  

 Let $\mathcal{X}$ be a real Banach space, and let 
${\bf T}=(T_1,T_2)\subset\mathcal{B(X)}$ be a pair of commuting
operators. The extended pair
${\bf T}_\C=(T_{1\C},T_{2\C})\subset\mathcal{B(X_\C)}$ also consists of commuting operators. For simplicity, we set 
$$
Q({\bf T}_\C):=\left(\begin{array}{cc} T_{1\C}, & T_{2\C}, \\ -T_{2\C}, & T_{1\C},
\end{array}\right)
$$
which acts on the complex Banach space $\mathcal{X}_\C^2$.

% Def 4

\begin{Def}\label{Q-jspectrum}\rm Let $\mathcal{X}$ be a real Banach space. For a given pair ${\bf T}=(T_1,T_2)\subset\mathcal{B(X)}$ of commuting operators, the set of those 
$Q(\z)\in\H,\,\z=(z_1,z_2)\in\C^2$, such that the operator
$$T_1^2+T_1^2-2\Re{z_1}T_1-2\Re{z_2}T_2+\vert z_1\vert^2+
\vert z_2\vert^2
$$
is invertible in $\mathcal{B(X)}$
is said to be the {\it quaternionic joint resolvent} (or simply the $Q$-{\it joint resolvent}) of ${\bf T}$. 

The complement $\sigma_\H({\bf T})=\H\setminus\rho_\H({\bf T})$ is called 
the {\it quaternionic joint spectrum} (or simply the $Q$-{\it joint spectrum}) of ${\bf T}$.
\end{Def}

% Lemma 10

\begin{Lem}\label{Q-jsp} A quaternion $Q(\z)$\,$(\z\in\C^2)$ is in the set $\rho_\H({\bf T})$
if and only if the operators $Q({\bf T}_\C)-Q(\z),\,Q({\bf T}_\C)-Q(\z^*)$ are invertible in $\mathcal{B}(\mathcal{X}_\C^2)$.
\end{Lem}

{\it Proof}\, The assertion follows from the equalities
$$
\left(\begin{array}{cc} T_{1\C}-z_1 & T_{2\C}-z_2 \\
-T_{2\C}+\bar{z}_2 & T_{1\C}-\bar{z}_1\end{array}\right)
\left(\begin{array}{cc} T_{1\C}-\bar{z}_1 & -T_{2\C}+z_2 \\
T_{2\C}-\bar{z}_2 & T_{1\C}-z_1\end{array}\right)=
$$
$$
\left(\begin{array}{cc} T_{1\C}-\bar{z}_1 & -T_{2\C}+z_2 \\
T_{2\C}-\bar{z}_2 & T_{1\C}-z_1\end{array}\right)
\left(\begin{array}{cc} T_{1\C}-z_1 & T_{2\C}-z_2 \\
-T_{2\C}+\bar{z}_2 & T_{1\C}-\bar{z}_1\end{array}\right)=
$$
$$
[(T_{1\C}-z_1)(T_{1\C}-\bar{z}_1)+ 
(T_{2\C}-z_2)(T_{2\C}-\bar{z}_2)]{\bf I}.
$$
for all $(\z=(z_1,z_2)\in\C^2)$.

Consequently, the operators $Q({\bf T}_\C)-Q(\z),\,Q({\bf T}_\C)-Q(\z^*)$ are invertible in
$\mathcal{B}({\mathcal X}_\C^2)$ if and only if the operator 
$(T_{1\C}-z_1)(T_{1\C}-\bar{z}_1)+ 
(T_{2\C}-z_2)(T_{2\C}-\bar{z}_2)$
is invertible in  $\mathcal{B}(\mathcal{X}_\C)$.  Because we have
$$
 T_{1\C}^2+T_{2\C}^2-2\Re{z_1}T_{1\C}-2\Re{z_2}T_{2\C}+\vert z_1\vert^2+\vert z_2\vert^2=
$$
$$ 
[T_1^2+T_1^2-2\Re{z_1}T_1-2\Re{z_2}T_2+\vert z_1\vert^2+
\vert z_2\vert^2]_\C,
$$
the operators $Q({\bf T}_\C)-Q(\z),\,Q({\bf T}_\C)-Q(\z^*)$ are invertible in
$\mathcal{B}({\mathcal X}_\C^2)$ if and only if the operator 
$T_1^2+T_1^2-2\Re{z_1}T_1-2\Re{z_2}T_2+\vert z_1\vert^2+
\vert z_2\vert^2$
is invertible in $\mathcal{B(X)}$.
\medskip

Lemma \ref{Q-jsp} shows that the set $\sigma_\H({\bf T})$ is 
$*$-{\it invariant}, that is $q\in\sigma_\H({\bf T})$ if and only if $q^*\in\sigma_\H({\bf T})$. Putting 
$$
\sigma_{\C^2}({\bf T}):=\{\z\in\C^2;Q(\z)\in\sigma_\H({\bf T}\},
$$ 
the set $\sigma_{\C^2}({\bf T})$ is also $*$-invariant, that is,
$\bf z\in\sigma_{\C^2}({\bf T})$ if and only if $\bf z^*\in\sigma_{\C^2}({\bf T})$

% Rem 20

\begin{Rem}\rm For the extended pair 
${\bf T}_\C=(T_{1\C},T_{2\C})\subset\mathcal{B(X_\C)}$ of the commuting pair ${\bf T}=(T_1,T_2)\subset\mathcal{B(X)}$ there is an 
interesting connexion with the {\it joint spectral theory} of 
J. L. Taylor (see \cite{Tay}; see also \cite{Vas3}). Namely, if the operator
$T_{1\C}^2+T_{2\C}^2-2\Re{z_1}T_{1\C}-2\Re{z_2}T_{2\C}+\vert z_1\vert^2+\vert z_2\vert^2$ is invertible, then the point
$\z=(z_1,z_2)$ belongs to the joint resolvent of ${\bf T}_\C$. 
Indeed, setting 
$$
R_j({\bf T}_\C,\z)=(T_{j\C}-\bar{z}_j)(T_{1\C}^2+T_{2\C}^2-2\Re{z_1}T_{1\C}-2\Re{z_2}T_{2\C}+\vert z_1\vert^2+\vert z_2\vert^2)^{-1},
$$ $q=Q(\z)$
for $j=1,2$, we clearly have
$$
(T_{1\C}-z_1)R_1({\bf T}_\C,\z)+(T_{2\C}-z_2)R_2({\bf T}_\C,\z)
={\bf I},
$$ 
which, according to \cite{Tay}, implies that $\z$ is in the 
joint resolvent of ${\bf T}_\C$. A similar argument shows that,
in this case the point  $\z^*$ also belongs to the joint resolvent of ${\bf T}_\C$. In addition, if $\sigma(T_\C)$
designates the Taylor spectrum of $T_\C$, we have the 
inclusion $\sigma(T_\C)\subset\sigma_{\C^2}({\bf T})$.                                                                                                                                                                                                                                                                                                                                                                                                                                                     
In particular, for every complex-valued function $f$ analytic in a  neighborhood of $\sigma_{\C^2}({\bf T})$, the operator 
$f(\bf T_\C)$ can be computed via Taylor's analytic functional
calculus. In fact, we have a Martinelli type formula for the
analytic functional calculus:
\end{Rem}

\begin{Thm} Let $\mathcal{X}$ be a real Banach space, let ${\bf T}=(T_1,T_2)\subset\mathcal{B(X)}$ be a pair of commuting operators, let $U\subset\C^2$ be an open set, let  
$D\subset U$ be a $*$-invariant bounded domain
containing $\sigma_{\C^2}({\bf T})$, with piecewise-smooth boundary $\Sigma$, and let $f\in\mathcal{O}(U)$. Then we have
$$
f({\bf T}_\C)=\frac{1}{(2\pi i)^2}\int_\Sigma f(\z))L({\bf z,T_\C})^{-2}(\bar{z}_1-T_{1\C})d\bar{z}_2-(\bar{z}_2-T_{2\C})
d\bar{z}_1]dz_1 dz_2,
$$
where 
$$
L({\bf z,T_\C})=T_{1\C}^2+T_{2\C}^2-2\Re{z_1}T_{1\C}-2\Re{z_2}T_{2\C}+\vert z_1\vert^2+\vert z_2\vert^2.
$$
\end{Thm}

{\it Proof.}\, Theorem III.9.9 from \cite{Vas3} implies that the 
map $\mathcal{O}(U)\ni f\mapsto f({\bf T}_\C)\in\mathcal{B(X_\C)}$, defined in terms of Taylor's analytic functional calculus, is unital, 
linear, multiplicative, and ordinary complex polynomials in $\z$
are transformed into polynomials in ${\bf T}_\C$ by simple substitution, where $\mathcal{O}(U)$ is the algebra of all analytic functions in the open set $U\subset\C^2$, provided 
$U\supset\sigma({\bf T}_\C)$.

The only thing to prove is that, when $U\supset\sigma_{\C^2}({\bf T})$, Taylor's functional calculus is given by the stated (canonical) formula. In order to do that, we use an argument  from the proof of Theorem 
III.8.1 in \cite{Vas3}, to make explicit the integral III(9.2) from \cite{Vas3}.

We consider the exterior algebra
$$
\Lambda[e_1,e_2,\bar{\xi_1},\bar{\xi_2},\mathcal{O}(U)\otimes\mathcal{X}_\C]=
\Lambda[e_1,e_2,\bar{\xi_1},\bar{\xi_2}]\otimes\mathcal{O}(U)\otimes\mathcal{X}_\C,
$$
where the indeterminates $e_1,e_2$ are to be associated with the pair ${\bf T}_\C$, we put $\bar{\xi_j}=d\bar{z}_j,\,j=1,2$, and consider the operators
$\delta=(z_1-T_{1\C})\otimes e_1+(z_2-T_{2\C})\otimes e_2,\,\bar{\partial} = 
(\partial/\partial\bar{z_1})\otimes\bar{\xi_1}+(\partial/\partial\bar{z_2})\otimes\bar{\xi_2}$, acting naturally on this 
exterior algebra, via the calculus with exterior forms. 

To simplify the computation, we omit the symbol $\otimes$, and the exterior product will be denoted simply par juxtaposition. 

We fix the exterior form $\eta=\eta_2=fye_1e_2$ for some 
$f\in\mathcal{O}(U)$ and  $y\in\mathcal{X}_\C$, which clearly satisfy the equation 
$(\delta+\bar{\partial})\eta=0$, and look for a solution
$\theta$ of the equation $(\delta+\bar{\partial})\theta=\eta$. 
We write $\theta=\theta_0+\theta_1$, where $\theta_0,\theta_1$
are of degree $0$ and $1$ in $e_1,e_2$, respectively. Then the equation $(\delta+\bar{\partial})\theta=\eta$ can be written 
under the form $\delta\theta_1=\eta,\,\delta\theta_0=-\bar{\partial}\theta_1$, and $\bar{\partial}\theta_0=0$. 
Note that
$$
\theta_1=fL({\bf z,T_\C})^{-1}[(\bar{z}_1-T_{1\C})ye_2-
(\bar{z}_2-T_{2\C})]ye_1
$$
is visibly a solution of the equation $\delta\theta_1=\eta$. 
Further, we have
$$
\bar{\partial}\theta_1=fL({\bf z,T_\C})^{-2}
[(z_1-T_{1\C})(\bar{z}_2-T_{2\C})y\bar{\xi}_1e_1-
(z_1-T_{1\C})(\bar{z}_1-T_{1\C})y\bar{\xi}_2e_1+
$$
$$
(z_2-T_{2\C})(\bar{z}_2-T_{2\C})y\bar{\xi}_1e_2-
(z_2-T_{2\C})(\bar{z}_1-T_{1\C})y\bar{\xi}_2e_2]=
$$
$$
\delta[fL({\bf z,T_\C})^{-2}(\bar{z}_1-T_{1\C})y\bar{\xi}_2-
fL({\bf z,T_\C})^{-2}(\bar{z}_2-T_{2\C})y\bar{\xi}_1],
$$
so we may define
$$
\theta_0=-fL({\bf z,T_\C})^{-2}(\bar{z}_1-T_{1\C})y\bar{\xi}_2+
fL({\bf z,T_\C})^{-2}(\bar{z}_2-T_{2\C})y\bar{\xi}_1.
$$
Formula III(8.5) from \cite{Vas3} shows that 
$$
f({\bf T}_\C)y=-\frac{1}{(2\pi i)^2}\int_U\bar{\partial}(\phi\theta_0)dz_1 dz_2=
$$
$$
\frac{1}{(2\pi i)^2}\int_\Sigma f(\z))L({\bf z,T_\C})^{-2}[(\bar{z}_1-T_{1\C})yd\bar{z}_2-(\bar{z}_2-T_{2\C})y
d\bar{z}_1]dz_1 dz_2,
$$
for all $y\in\mathcal{X}_\C$, via Stokes's formula, where $\phi$ is a smooth function such 
that $\phi=0$ in a neighborhood of $\sigma_{\C^2}({\bf T})$,
$\phi=1$ on $\Sigma$ and the support of $1-\phi$ is compact.

\begin{Rem}\rm (1) We may extend the previous functional calculus to $\mathcal{B(X}_\C)$-valued analytic functions, setting, for such a function $F$ and with the notation from above,
$$
F({\bf T}_\C)=\frac{1}{(2\pi i)^2}\int_\Sigma F(\z))L({\bf z,T_\C})^{-2}(\bar{z}_1-T_{1\C})d\bar{z}_2-(\bar{z}_2-T_{2\C})
d\bar{z}_1]dz_1 dz_2.
$$
In particular, if $F(\z)=\sum_{j,k\ge0}A_{jk\C}z_1^jz_2^k$,
with $A_{j,k}\in\mathcal{B(X)}$, where the series is convergent in neighborhood of $\sigma_{\C^2}({\bf T})$, we may define 

$$F({\bf T}):=F({\bf T}_\C)\vert\mathcal{X}=\sum_{j,k\ge0}A_{jk}T_1^jT_2^k\in\mathcal{B(X)}.$$
 
(2) The connexion of the spectral theory of pairs with the algebra of quaternions is even stronger in the case of complex Hilbert
spaces. Specifically, if $\mathcal{H}$ is a complex Hilbert space and ${\bf V}=(V_1,V_2)$ is a commuting pair of bounded
linear operators on $\mathcal{H}$, a point $\z=(z_1,z_2)\in\C^2$
is in the joint resolvent of ${\bf V}$ if and only if the 
operator $Q({\bf V})-Q(\z)$ is invertible in $\mathcal{H}^2$,
where 
$$
Q({\bf V})=\left(\begin{array}{cc} V_1 & V_2 \\ -V_2^* & V_1^*
\end{array}\right).
$$
(see \cite{Vas1} for details). In this case, there is also a Martinelli type formula which can be used to construct the associated analytic functional calculus (see \cite{Vas2},\cite{Vas3}). An approach to such a construction in Banach spaces, by using a so-called splitting joint spectrum, can be 
found in \cite{MuKo}.  
\end{Rem}

\end{document}